\documentclass[11pt]{amsart}
\usepackage{amsthm,amsfonts,amssymb,amsmath,oldgerm}
\numberwithin{equation}{section}
\numberwithin{figure}{section}
\numberwithin{table}{section}
\usepackage{fullpage}
\usepackage{amsmath}
\usepackage{setspace}
\usepackage{stmaryrd}
\usepackage{mathrsfs}
\usepackage{fancyhdr}
\usepackage{graphicx}
\usepackage{enumerate}
\usepackage{bm}
\usepackage{bbm}
\usepackage{dsfont}
\usepackage{multirow}
\usepackage{psfrag}
\usepackage[font=scriptsize]{caption}
\usepackage[font=scriptsize]{subcaption}
\usepackage{listings}
\usepackage{tikz}
\usepackage{empheq}
\usepackage{cases}
\usetikzlibrary{matrix} 

\usepackage{hyperref}
\hypersetup{
    colorlinks,               % color of internal links (change box color with linkbordercolor)
    citecolor=blue,        % color of links to bibliography
    filecolor=blue,         % color of file links
    urlcolor=blue,         % color of external links
    linkcolor=blue,
  }

\newcommand\dD{\mathrm{d}}
\newcommand\ds{\displaystyle}
\newcommand{\N}{{\mathbb N}}
\newcommand{\C}{{\mathbb C}}
\newcommand{\R}{{\mathbb R}}

\newcommand\cC{\mathcal C}
\newcommand\cD{\mathcal D}
\newcommand\cE{\mathcal E}

\newcommand\cJ{\mathcal J}

\newcommand\cK{\mathcal K}
\newcommand\cR{\mathcal R}
\newcommand\cP{\mathcal P}
\newcommand\cS{\mathcal S}
\newcommand\cF{\mathcal F}

\newcommand{\mC}{{\mathscr C}}
\newcommand{\mP}{{\mathscr P}}
\newcommand{\be}{\begin{equation}}
\newcommand{\ee}{\end{equation}}

\newtheorem{theorem}{Theorem}[section]

\newtheorem{remark}[theorem]{Remark}


\begin{document}

\title[On the approximation of the von-Neumann equation]{On the
  approximation of the von-Neumann equation in the semi-classical
  limit. Part I : numerical algorithm}

\author[Francis Filbet]{Francis Filbet}
\address[F.F.]{IMT, Université Paul Sabatier, 31062 Toulouse Cedex, France}
\email{francis.filbet@math.univ-toulouse.fr}

\author[François Golse]{Fran\c cois Golse}
\address[F.G.]{CMLS, \'Ecole polytechnique, 91128 Palaiseau Cedex, France}
\email{francois.golse@polytechnique.edu}

\begin{abstract}
  We propose a new approach to discretize the von Neumann equation, which is efficient in the semi-classical limit. This method
  is first based on the so called Weyl's variables to address the
  stiffness associated with the equation. Then, by applying a
  truncated Hermite expansion of the density operator, we successfully
  handle this stiffness. Additionally, we
  develop a finite volume approximation for practical
  implementation and conduct numerical simulations to illustrate the
  benefits of our approach. This asymptotic preserving
  numerical approximation, combined with the use of Hermite
  polynomials, provides a useful tool  for solving
  the von Neumann equation in all regimes, near classical or not.
  \end{abstract}

\date{\today}

\maketitle


\tableofcontents

%%%%%%%%%%%%%%%%%%%%%%%%%%%%%%%%%%%%%%%%%%%%%%%%%%%%%%%%%%%%%%%%%%%%%%%%%%%

\section{Introduction}
\label{sec:1}
\setcounter{equation}{0}
\setcounter{figure}{0}
\setcounter{table}{0}

%%%%%%%%%%%%%%%%%%%%%%%%%%%%%%%%%%%%%%%%%%%%%%%%%%%%%%%%%%%%%%%%%%%%%%%%%%%

We consider an operator $\hat\rho$, called the density-matrix
operator, satisfying the Liouville-von Neumann or simply von Neumann equation in the operator formulation
\begin{equation*}
\left\{
\begin{array}{l}
\ds i\,\hbar \,\partial_t \hat\rho \,=\, [H,\hat\rho],
\\[1.1em]
\hat\rho(0)\,=\,\hat\rho^{in}\,,
\end{array}\right.
\end{equation*}
where $\hbar$ is the reduced Planck constant, the operator $\hat\rho$
{is acting on the  Hilbert space $L^2(\R^d)
:=L^2(\R^d,\C)$ taken with the inner product
$$
\langle f, g\rangle = \int_{\R^d} f(x)\,\overline{g(x)}\,\dD x,
$$}
and the associated  norm of a state to be $\|.\| =
\sqrt{\langle.,.\rangle}$. Moreover, the quantum mechanical Hamiltonian $H$ is given, for instance, by
$$
H = -\frac{\hbar^2}{2m}\,\Delta  \,+\, V (X),
$$
where $V$ a real-valued function such that $H$ is a self-adjoint operator on $L^2(\R^d)$. On the other hand, 
for any $\phi\in L^2(\R^d)$, $\phi \mapsto V \phi(X) \,=\, V (X)\,\phi(X)$
while $[H,\hat\rho] = H \,\hat\rho - \hat\rho\,H$ is the
commutator. Here,  we suppose that  $\hat\rho(t)$ is a density operator on $L^2(\R^d)$, such that 
\begin{equation}
\label{prop:00}
\hat\rho(t)^*\,=\,\hat\rho(t)\,\ge\,
0\,,\quad\text{Tr}(\hat\rho(t))\,=\,1, \quad t \geq 0.
\end{equation}
In particular, $\hat \rho(t)$ is trace-class, and therefore
Hilbert-Schmidt, hence  by \cite[Theorem 6.12]{brezis},  there exists a
unique $\rho(t|\,X,Y)\in {L^2_{X,Y}:= L^2(\R^d\times \R^d, \C)}$, such that for any $\phi\in L^2(\R^d)$, we
define $\hat\rho(t)\phi$ as the function
$$
X\,\mapsto \, \hat\rho(t)\,\phi(X)\,=\;\int_{\R^d}\rho(t|\,X,Y)\,\phi(Y)\,\dD Y\,.
$$
Thus, the von Neumann equation can be written as follows for the
function $(X,Y)\mapsto \rho(t |\,X,Y)$,
\begin{equation}
  \label{vNn}
  \left\{
\begin{array}{l}
\ds i\,\hbar\,\partial_t \rho(t|\,X,Y) \,=\,-\frac{\hbar^2}{2m}\,\left(\Delta_X\,-\,\Delta_Y\right)\,\rho(t|\,X,Y)\,+\,\left(V(X)\,-\,V(Y)\right)\,\rho(t|\,X,Y)\,,
\\[1.1em]
 \ds     \rho(0)\,=\,\rho^{in}\,.
    \end{array}\right.
\end{equation}

The statement worth noting is that for $\hbar \ll 1$, the dynamics
described by the von Neumann equation becomes stiff, leading to a
multiscale problem. Thus a major challenge in
quantum dynamics  is to develop efficient
numerical techniques capable of handling a large variety of wavelengths. For instance, numerical  schemes used to solve the
Schr\"{o}dinger equation typically demand that both the time step
$\Delta t$ and mesh size, in the semiclassical regime when
$\hbar\ll 1$, are of order $O(\hbar)$, or may even be smaller than
$\hbar$ \cite{markowich1999, athanass}.  Actually, time splitting methods have
the advantage of using  large time steps  when only the physical
observables are concerned \cite{BaoJinMarko}. We also refer to \cite{russo2013,russo2014,  jin2022,miao2023} and the references therein for the
numerical approximation of the Schr\"{o}dinger equation in the
semi-classical limit.  A fundamental concept to
understand these mesh strategies is the Wigner transform, which serves
as a useful tool for investigating the semiclassical limit of the
Schr\"{o}dinger equation.  Additionally for the von Neumann equation, the density operator, which is Hermitian,
positive semi-definite, and has  trace  equal to one by definition, should
ideally maintain these properties at the discrete level in order to ensure realistic
outcomes.

In this context, a first approach  has been proposed in \cite{hellsing1986efficient} and consists of successive application
of short time propagators which are evaluated by using fast Fourier
transforms. This method demonstrates high efficiency  for short time when the density matrix is
spatially localized but, unfortunately, it requires a small  time step
of order $\hbar$.   Later in
\cite{berman1992solution}, the authors proposed a numerical scheme
for a related model, again based on the Fourier pseudo-spectral
method. It allows a description both in configuration as well as in
momentum space. More recently, various structure preserving  schemes
have been studied to guarantee  trace conservation and positivity of
the discrete density matrix,  allowing numerical simulations over
large time intervals.

Actually, as we shall see more clearly later,  the von Neumann equation is strongly related to the Wigner
equation, which shares the same kind of numerical difficulties arising from the nonlocal and highly-oscillating pseudo-differential
operator. It has  been shown that the Wigner function can be
discretized in a robust way using adaptive pseudo-spectral methods \cite{shao2011adaptive,xiong2016advective}. In these methods, the
carefully crafted spectral decomposition of the Wigner function
enables the oscillatory components introduced by the Wigner kernel to
be solved exactly. Finally, the oscillatory quantum effects can also
be mitigated by decomposing the potential in a classical and quantum
part \cite{sellier2015comparison, sels2013wigner} or by reformulating
the Wigner equation using spectral components of the force field \cite{jcp2017}.

Our aim is to design numerical schemes for  the von Neumann equation \eqref{vNn}  allowing an evaluation of
the main observables (mass, momentum, and energy densities...) without
having to use time steps or space discretization in  $O(\hbar)$. We already know that it is possible to choose  time steps much larger
than $\hbar$ and get error estimates uniform  in  $\hbar$
\cite{BaoJinMarko,FGJinPaul} by applying time splitting schemes. But
the discussion in \cite{BaoJinMarko,FGJinPaul} leaves aside the space
discretization issue, which is the core problem addressed here.

Our first move is to write the von Neumann equation in terms of  Weyl's
variables (see \eqref{WeylVar}). This removes all stiffness as $\hbar\ll 1$, and the new
equation so obtained is well adapted to the semi-classical
limit as $\hbar$ tends to zero. Then, in Section \ref{sec:2}, we expand the 
density operator using Hermite polynomials leading to a spectrally accurate
approximation. We shall briefly review  error estimates, proven in our
forthcoming work \cite{FG2:24}, which are uniform in the parameter $\hbar$.  Next, in Section \ref{sec:3} we propose a
finite volume scheme for the space discretization coupled with a Crank-Nicolson time
discretization. Furthermore, in order to reduce the computational
effort when $\hbar \ll 1$, we present an approximation of the
non-local operator obtained after our change of variables. Finally, we
perform extensive numerical simulations with various potentials in
order to illustrate the efficiency of our numerical scheme (see Section \ref{sec:4}).   

%%%%%%%%%%%%%%%%%%%%%%%%%%%%%%%%%%%%%%%%%%%%
% FF : completer par le plan du papier
%%%%%%%%%%%%%%%%%%%%%%%%%%%%%%%%%%%%%%%%%%%%

\subsection{Quantum dynamics in  Weyl's variables}
\label{sec:1.1}
%%%%%%%%%%%%%%%%%%%%%%%%%%%%%%%%%%%%%%%%%%%%%%
%
%              Quantum dynamics in  Weyl's variable       
%                       
%%%%%%%%%%%%%%%%%%%%%%%%%%%%%%%%%%%%%%%%%%%%%%

To remove the stiffness due to the parameter $\hbar$ in the von
Neumann equation, the key idea consists in rewriting
\eqref{vNn} in the new variables defined below,  called “Weyl's variables”:
\begin{equation}
  \label{WeylVar}
  \left\{
    \begin{array}{l}
      \ds x\,:=\,\frac{X\,+\,Y}{2}\,,\\[0.9em]
      \ds y\,:=\,\frac{X\,-\,Y}{\hbar}\,,
    \end{array}\right.
\end{equation}
which can be explicitly inverted as
$$
X\,=\,x\,+\,\frac{\hbar}{2}\,y\,,\qquad Y\,=\,x\,-\,\frac{\hbar}{2}\,y\,.
$$
Although this terminology is not common, it is obviously suggested by
Weyl's  quantization (see formula (1.1.9) of
\cite{Lerner}). Hence, we set
$$
R(t|\,x,y)\,:=\, \rho\left(t | \,X,Y\right)\,;
$$
and find that $\rho$ is solution to \eqref{vNn} if and only if
$R$ satisfies
$$
 \left\{
     \begin{array}{l}
\ds i\,\hbar\,\partial_t R(t|\,x,y)\,=\,-\hbar\,\sum_{j=1}^d\partial_{x_j}\partial_{y_j}
R(t|\,x,y)\,+\,\left(V\left(x+\frac\hbar 2\, y\right)\,-\,V\left(x-\frac\hbar 2\, y\right)\right)\,R(t|\,x,y)\,,
\\[1.1em]
   \ds   R(0)\,=\,R^{in}_\hbar\,,
     \end{array}\right.
   $$
or  equivalently
\begin{equation}
  \label{vNvW}
   \left\{
     \begin{array}{l}
       
\ds \partial_t R(t|\,x,y)\,=\,i\,\sum_{j=1}^d\partial_{x_j}\partial_{y_j}
R(t|\,x,y)\,-\,i\frac{V\left(x+\frac\hbar 2\,
       y\right)\,-\,V\left(x-\frac\hbar 2\, y\right)}{\hbar}\,
       R(t|\,x,y)\,,
\\[1.1em]
      R(0)\,=\,R^{in}_\hbar\,.
      \end{array}\right.
  \end{equation}
  
Let us first review some elementary properties on the operator $\hat\rho$ and
their interpretation in terms of the function $R$. We first recall that $\hat\rho$ is a density operator on $L^2(\R^{2d})$, that is
$\hat\rho (t)^*=\hat\rho (t)$, which means that  $R$ satisfies 
\begin{equation}
\label{SAvW}
R(t|\,x,-y)\,=\,\overline{R(t|\,x,y)}\,,
\end{equation}
whereas the trace property becomes 
\begin{equation}
 \label{TRvW}
\text{Tr}\left(\hat\rho (t)\right) \,=\,\int_{\R^d}\rho(t|\,X,X)\,\dD X\,=\,\int_{\R^d}
R(t|\,x,0)\,\dD x\,=\,1\,.
\end{equation}
Translating the positivity condition on $\hat\rho(t)\geq 0$ in terms
of the function $R$ is more difficult and we leave it aside (see however in \cite{GolseMoller} a characterization of rank-$1$ density operators in Weyl's variables).

On the one hand, {since $\hat\rho(t)=\exp(-itH/\hbar)\hat\rho(0)\exp(itH/\hbar)$ for all $t$, one has
\[
\hat\rho(t)^*=\hat\rho(t)\text{ since }\hat\rho(0)^*=\hat\rho(0)\quad\text{ and }\mathrm{Tr}(\hat\rho(t)^2)=\mathrm{Tr}(\hat\rho(0)^2)<\infty\,,
\]
and it is equally well-known that the Hilbert-Schmidt norm of an operator is the $L^2$ norm of its integral kernel, so that
\[
\|\hat\rho(t|\cdot,\cdot)\|^2_{L^2(\mathbf R^d\times\mathbf R^d)}=\mathrm{Tr}(\hat\rho(t)^2)=\mathrm{Tr}(\hat\rho(0)^2)=\|\hat\rho(0|\cdot,\cdot)\|^2_{L^2(\mathbf R^d\times\mathbf R^d)}\,.
\]
The simple change of variables \eqref{WeylVar}, leading to
\[
\|\hat\rho(t|\cdot,\cdot)\|^2_{L^2(\R^d\times\R^d)}\,=\,\hbar^d\| R(t|\cdot,\cdot)\|^2_{L^2(\R^d\times\R^d)}\;,
\]}
which yields to  the conservation of the $L^2$ norm for
equation \eqref{vNvW} : for each $t\in[0,+\infty)$,
$$
%\label{esti:L2}
  \iint_{\R^{2d}}|R(t|\,x,y)|^2\,\dD x\,\dD y\,=\,\| 
R^{in}_\hbar\|^2_{L^2}\,<\,\infty\,.
$$
This latter property is crucial in the stability analysis  of a
numerical method \cite{FG2:24}.

%%%%%%%%%%%%%%%%%%%%%%%%%%%%%%%%%%%%%%%%%%%%%%%%%%%%%%%%%%%%%%%%%%%%%%%%%%%%%%%%%%%%%%%%%%%%%%
%
%    WIGNER transform
%
%%%%%%%%%%%%%%%%%%%%%%%%%%%%%%%%%%%%%%%%%%%%%%%%%%%%%%%%%%%%%%%%%%%%%%%%%%%%%%%%%%%%%%%%%%%%%%

Finally, let us point out that the von Neumann equation written in Weyl's variables \eqref{vNvW} is well
adapted to making the link with other formulations, such as  the Wigner
equation. Indeed,  for all  $(t,x)\in\R^+\times\R^d$, we  define  the Wigner function
$W(t,x,\cdot)$ by taking  the Fourier
transform of the function $R(t|\,x,\cdot)$
$$
W(t,x,\xi) \,=\, \frac{1}{(2\pi)^d}\int_{\R^d} R(t|\,x,y) \, e^{-i\,y\cdot
  \xi}\,\dD y\,. 
$$

The integral on the right-hand side is  not absolutely converging in general 
({\it i.e.} not a Lebesgue integral), but is the  Fourier(-Plancherel) transform of an $L^2$ function on $\R^d$.

Then the function $W$ is formally a solution of the so-called Wigner
equation \cite[Chapter 11]{jungel}
\begin{equation}
\label{eq:Wigner}
  \partial_t W \,+\,\xi\cdot\nabla_x W \,-\, \Theta[V]\, W\,=\,0\,,
\end{equation}
where the operator $\Theta[V]$ is given by
{$$
\Theta[V] W(t,x,\xi) \,=\,  \,\frac{i}{(2\pi)^d\,\hbar} \,\iint_{\R^{2d}}\left( V\left(x+\frac\hbar 2
  y\right)\,-\,V\left(x-\frac\hbar2 y\right) \right) \,W(t,x,\xi^\prime)\,
e^{-i\,y\cdot (\xi-\xi^\prime)}\,\dD\xi^\prime\,\dD y\,. 
$$}
{It is worth mentioning that while $R$ is a complex function, it
satisfies $R(t|\,x, -y) = \overline{R(t|\,x, y)}$, so that its Fourier
transform in $y$,  $W(t,x,.)$ is real-valued \cite[pages
564]{LionsPaul} or \cite{GMMP,AMP:2009}.}
The local term $\xi\cdot\nabla_x W$  is the quantum analogue of the
classical transport term of the Liouville equation whereas  the
nonlocal term $\Theta[V ]\,W $ describes the influence of the electric
potential. In other words, while the classical limit of the von Neumann equation
is a singular perturbation problem in the original variables, it
becomes a regular perturbation problem for smooth potentials, when
written in Weyl's variables, or for the Wigner equation.

\subsection{Semi-classical limit}
\label{sec:1.2}
The formulation of the von Neumann
equation using Weyl's variables facilitates the semi-classical
limit. Indeed,  passing to the limit  as $\hbar\to 0$ in the latter equation  does not lead
to any difficulty related to the stiffness in $\hbar$. Hence,  since $V\in
\mC^1(\R^d)$, we have
$$
\frac{V(x+\frac\hbar 2\, y)-V(x-\frac\hbar 2\, y)}{\hbar}\to\nabla V(x)\cdot y\qquad\text{as }\hbar\to 0\,.
$$
Therefore, we set
\begin{equation}
  \label{def:Eh}
\cE^\hbar(x,y) \,=\,  \frac{V(x+\frac\hbar 2\, y)-V(x-\frac\hbar 2
\,  y)}{\hbar}\,-\, \nabla V(x)\cdot y\,,
\end{equation}
and $R$ solves
\begin{equation}
  \label{vNvW2}
   \left\{
    \begin{array}{l}
\ds\partial_t R(t|\,x,y)\,=\,i\,\sum_{j=1}^d\partial_{x_j}\partial_{y_j}
R(t|\,x,y)\,-\,i\left(\nabla V(x)\cdot y  \,+\,\cE^\hbar(x,y)\right)
\,R(t|\,x,y)\,,
\\[1.1em]
      R(0)\,=\,R^{in}_\hbar\,,
      \end{array}\right.
\end{equation}
where $\cE^\hbar$ tends to zero as $\hbar\rightarrow 0$. {This latter
equation is equivalent to \eqref{vNvW} and it will the starting point
of the design of our numerical scheme.} 

The solution obtained
by neglecting the source term $\cE^\hbar$ corresponds to the
semi-classical limit, given by
\begin{equation}
  \label{eq:lim}
  \left\{
    \begin{array}{l}
\ds\partial_t\hat
R(t|\,x,y)\,=\,i\,\sum_{j=1}^d\partial_{x_j}\partial_{y_j}\hat
R(t|\,x,y)\,-\,i \,\nabla V(x)\cdot y  \,\hat R(t|\,x,y)\,,
\\[1.1em]
      \hat R(0)\,=\,\hat R^{in}_\hbar\,.
      \end{array}\right.
  \end{equation}
{The error made
in approximating the solution $R(t|\,x, y)$ of
\eqref{vNvW2}   by  the solution $\hat
R(t|\,x, y)$ of \eqref{eq:lim} }  is analyzed in Theorem \ref{th:01} below.
Moreover, equation \eqref{eq:lim} can be interpreted using  the Wigner function
$\widehat W(t,x,.)$ computed from $\hat R(t|\,x,\cdot)$
$$
\widehat W(t,x,\xi) \,=\, \frac{1}{(2\pi)^d}\,\int_{\R^d} \hat R(t|\,x,\xi) \, e^{-i\,y\cdot
  \xi}\,\dD y\,.
$$
Taking the Fourier transform with respect to  the $y$ variable  in
\eqref{eq:lim},  we see that the function $\widehat W$ solves the Liouville equation 
\begin{equation*}
  \partial_t \widehat W \,+\,\xi\cdot\nabla_x \widehat W \,-\, \nabla_x V\cdot
  \nabla_\xi \widehat W\,=\,0\,.
\end{equation*}
 This limit was made rigorous by P.-L. Lions and T. Paul
 \cite{LionsPaul}, and  P. Markowich and C. Ringhofer \cite{MN89,MR89}
for smooth potentials. {We also mention \cite{AP:2011}  for strong convergence}. The limit was also performed in \cite[Section 1.4]{ref2}
by an asymptotic expansion of $\cE^\hbar$ in powers of $\hbar$. For references
on the mathematical analysis of the Wigner equation (or the Wigner–Poisson system), we refer to the review of A. Arnold \cite{AA}.

In the following, we restrict ourselves to the one-dimensional case
$d=1$. Let us fix an integer $m\geq 0$ and suppose that the potential $V$ is such that
\begin{equation}
  \label{hyp:V1}
  V\in \mC^{m+1}(\R) % \quad{\rm and
                                % }\quad V(x)\to+\infty, \quad {\rm as}\quad |x|\to\infty\,,
\end{equation}
and for all $n=2,\ldots,m+1$, there exists a constant $C>0$, not
depending on $\hbar$, such that 
\begin{equation}
  \label{hyp:V2}
\| V^{(n)}\|_{L^\infty}  \leq C,{\quad {\rm with\, } V^{(n)} =
\frac{d^n V}{d x^n}}. 
\end{equation}

From this assumption we get our first result on error estimates with
respect to $\hbar$ between the solution $R$ of the von Neumann
equation \eqref{vNvW} and $\hat R$ of its classical limit \eqref{eq:lim}.

\begin{theorem}
  \label{th:01}
 Assume that $V$ satisfies \eqref{hyp:V1}-\eqref{hyp:V2} with $m=3$ while  for all integers  $a$, $b$,
$\alpha$ and $\beta$ such that  $a+b+\alpha+\beta\le 3$, we suppose
that both initial data $R^{in}_\hbar$ and $\hat R^{in}_\hbar$ are such that
\begin{equation}
  \label{hyp:Rhat3}
  \left\{
    \begin{array}{l}
\ds\| x^a\partial^b_xy^\alpha\partial^\beta_y \hat R^{in}_\hbar \|_{L^2}  \,<
\, \cC_0\,,
      \\
      \ds\|  R^{in}_\hbar \|_{L^2}  \,<
      \, \cC_0\,,
      \end{array}\right.
    \ds \end{equation}
  where $\cC_0>0$ does not depend on $\hbar$ and $L^2=L^2_{xy}:=L^2(\R^2)$. Then  the solution $R$ to  the von Neumann equation \eqref{vNvW} and
the solution $\hat R$ to \eqref{eq:lim}  with the initial data
$R^{in}_\hbar$ and $\hat R^{in}_\hbar$ satisfy 
$$
\|R(t)\,-\,\hat R(t)\|_{L^2} \,\le\, \|R^{in}_\hbar\,-\,\hat R^{in}_\hbar\|_{L^2}
\,+\,\cC_t\,\hbar^2\,\qquad t\geq 0\,,
$$
where $\cC_t$ is a constant independent of $\hbar$.
\end{theorem}

The proof of this theorem is presented in \cite{FG2:24} and is a
simple  consequence of the regularity of the solution to the
semi-classical limit problem \eqref{eq:lim}.

\begin{remark}
Let us make the following comments on Theorem \ref{th:01}:
  \begin{itemize}
\item  of course, assumptions \eqref{hyp:Rhat3}   are not the most general setting for
 the semi-classical limit, in particular they exclude pure WKB
 states;

 \item  from \eqref{hyp:Rhat3}, we know that
$(\hat R^{in}_\hbar)_{\hbar>0}$ is strongly compact in $L^2$, but
allows fast oscillations in $\rho(t|\,X,Y)$ in the original
variables. Indeed, for a given  $(\tilde x,\tilde p)\in\R^2$, let us remind that if $(t,X)\mapsto
\Psi(t,X,\tilde x,\tilde p)$ is a solution to the Schr\"{o}dinger
equation,
$$
i\,\hbar \,\partial_t\Psi \,=\, H \Psi\,,
$$
then,  for a given weight function $K$, we may define the density
matrix with integral kernel
$$
\rho(t|\,X,Y) = \int_{\R^2} K(\tilde x,\tilde p)\, \Psi(t,X,\tilde
x,\tilde p)\,\overline\Psi(t,Y,\tilde
x,\tilde p)\,\dD\tilde x\,\dD\tilde p\,,
$$
which is a solution to the von Neumann equation \eqref{vNn}. Therefore
in our framework,  we may consider mixed sates obtained by first choosing the WKB initial
condition centered at $\tilde x\in \R$ and with an initial energy  $H(\tilde x,\tilde p)$, 
$$
\Psi^{in}(X, \tilde x, \tilde p) \,=\, \frac{1}{\sqrt{2\pi}} \,
\exp\left(- \frac{|X-\tilde x|^2}{4} \right) \, \exp\left(-i\,\frac{\tilde
p\,X}{\hbar}\right)
$$
and a weight function $K$ given by
$$
K(\tilde x,\tilde p) \,=\, \frac{1}{\sqrt{s_\hbar}}\,\exp\left(-\frac{|\tilde p-p_0|^2}{2
    s_\hbar}\right)\, \delta(\tilde x-x_0),
$$
where $s_\hbar = 1-\hbar^2/4$. Then, we construct the density function
$\rho^{in}$ as
\begin{eqnarray*}
\rho^{in}(X,Y) &=& \int_{\R ^2} K(\tilde x,\tilde p)
\,\Psi^{in}(X,\tilde x,\tilde p)\,\overline{\Psi^{in}}(Y,\tilde x,\tilde p)
  \,\dD\tilde x\,\dD\tilde p\,,
  \\
  &=&  \frac{1}{\sqrt{2\pi}}\, \exp\left(- \frac{\left(\tfrac{1}{2} (X+Y) - x_0\right)^2}{2}
      -  \frac{(X- Y)^2}{2\hbar^2}
      -  \frac{i\,p_0 \,(X-Y)}{\hbar} \right)\,,
\end{eqnarray*}
that is,
$$
R^{in}(x,y) \,=\,   \frac{1}{\sqrt{2\pi}}\, \exp\left(- \frac{|
    x- x_0|^2+ y^2}{2}-  i\,p_0 y\right)\,.
$$
In Section \ref{sec:4}, we will consider such initial conditions.
\end{itemize}
\end{remark}

%%%%%%%%%%%%%%%%%%%%%%%%%%%%%%%%%%%%%%%%%%%%%%%%%%%%%%%%%%%%%%%%%%%%%%%%%%%
\section{Hermite spectral method}
\label{sec:2}
\setcounter{equation}{0}
\setcounter{figure}{0}
\setcounter{table}{0}

%%%%%%%%%%%%%%%%%%%%%%%%%%%%%%%%%%%%%%%%%%%%%%%%%%%%%%%%%%%%%%%%%%%%%%%%%%%

The purpose of this section is to present a formulation of the von
Neumann equation {\eqref{vNvW2}} written in Weyl's variables based on Hermite
polynomials. We first use Hermite polynomials in the $y$
variable and write the von Neumann  equation  \eqref{vNvW}  as an infinite hyperbolic system for the Hermite coefficients
depending only on time and on the space variable $x$. The idea is to apply a Galerkin
method keeping only a small finite set of orthogonal
polynomials rather than discretizing the operator $R$ in $y$ on a grid. The merit of using
an orthogonal basis like the so-called scaled Hermite basis has been
shown in \cite{Holloway1996, Schumer1998} and more recently in \cite{Filbet2020,ref:5,bessemoulin2022cv, bf1-2024,bf2-2024,b-2024} for the Vlasov-Poisson system.
We define a weight function $\omega$ as 
\begin{equation*}
\omega(y):=\pi^{-1/4}\,e^{-y^2/2}\,,\qquad y\in\R\,,
\end{equation*}
and the sequence of Hermite polynomials (called the ``physicist's Hermite
polynomials''),
\begin{equation}
  \label{def:Hk}
H_k(y)\,=\,(-1)^k\,\omega^{-2}(y)\,\partial_y^k\omega^2(y)\,,\qquad k\ge 0\,.
\end{equation}

In this context the family of Hermite functions  $\left(\Phi_k\right)_{k\in\N}$ defined as
\begin{equation*}
%  \label{A-fHermite}
\Phi_k(y)\,=\,\frac1{\sqrt{2^kk!}}\,\omega(y)\, H_k(y)\,,\qquad k\ge 0
\end{equation*}
is a complete,  orthonormal system for the classical $L^2$ inner product, that is,
\[
\int_{\R}\,
\Phi_{k}(y)\,\Phi_{l}(y)\,\dD y
\,=\,
\delta_{k,l}\,, \quad{\rm and}\quad \overline{{\rm Span}\{\Phi_k, \quad k\geq 0\}}= L^2(\R).
\]
Moreover, the sequence $
\ds
\left(
H_{k}
\right)_{k\in\N} 
$ defined in \eqref{def:Hk} satisfies the recursion relation,
\begin{equation*}
   \left\{
    \begin{array}{l}
\ds H_{-1}=0, \quad H_{0}=1,
      \\[0.9em]
 \ds H_{k+1}(y)\,=\, 2\,y\,H_k(y) \,-\,  2\,k\,H_{k-1}(y)\,,\quad k\,\geq\,0\,.
\end{array}
      \right.
      \end{equation*}
Let us also point out that the Hermite functions also verify the following relations
\begin{equation}
 \label{2-bis}
\Phi_k'(y)\,=\, -\sqrt{\frac{k+1}{2}}\,\Phi_{k+1}(y) \,+\, \sqrt{\frac{k}{2}}\,\Phi_{k-1}(y) \,,\quad\forall\, k\,\geq\,0\,,
\end{equation}
and
\begin{equation}
  \label{yPhi}
y\,\Phi_k(y)\,=\, \sqrt{\frac{k+1}{2}}\,\Phi_{k+1}(y) \,+\, \sqrt{\frac{k}{2}}\,\Phi_{k-1}(y) \,,\quad\forall\, k\,\geq\,0\,.
\end{equation}

We now set for any $N\geq 0$,
$$
\mP_{N}(\R) \,=\, {\rm Span}\{ \Phi_k, \quad k=0,\ldots,N\}\,, 
$$
and expand the solution to the von Neumann and 
semi-classical limit equations  \eqref{vNvW} and \eqref{eq:lim}  in terms of  Hermite functions in the
variable $y$.

\subsection{Hermite approximation of  the von Neumann equation}
\label{sec:2.1}
 
We consider the decomposition of $R$, solution to 
 \eqref{vNvW}, into its components 
$
R\,=\,(R_{k})_{k\in\N}
$
in the Hermite basis:
\begin{equation*}
R\left(t|\,x,y\right)
\,=\,
\sum_{k\in\N}\,
R_{k}
\left(t|\,x
\right)\,\Phi_{k}(y)\,,
\end{equation*}
where
$$
R_k(t|\,x) \,:=\,\int_{\R}R(t|\,x,y)\,\Phi_k(y)\,\dD y\,,\quad k\geq 0\,.
$$
Using the recursion relations \eqref{2-bis}-\eqref{yPhi},  for all $k\geq 0$,
\begin{eqnarray*}
\int_{\R}\partial_yR(t|\,x,y)\,\Phi_k(y)\,\dD
  y&=&-\int_{\R}R(t|\,x,y)\,\partial_y\Phi_k(y)\,\dD y
\\[1.1em]
&=&\ds -\sqrt{\frac{k}2}\,R_{k-1}(t|\,x)\,+\,\sqrt{\frac{k+1}2}\,R_{k+1}(t|\,x)\,,
\end{eqnarray*}
with $R _{-1}=0$, whereas
\begin{eqnarray*}
\int_{\R}R(t|\,x,y)\,y\,\Phi_k(y)\,\dD
  y&=&\ds \sqrt{\frac{k}2}\,
       R_{k-1}(t|\,x)\,+\,\sqrt{\frac{k+1}2}\,R_{k+1}(t|\,x)\,.
\end{eqnarray*}
Similarly,
$$
\int_{\R}\left\{\frac{V(x+\frac\hbar 2\, y)-V(x-\frac\hbar 2\, y)}{\hbar}
  - V^\prime(x)\,y\right\}\,
R(t|\,x,y)\Phi_k(y)\,\dD y\,=\,\sum_{l\in\N}\cE^\hbar_{k,l}(x)\,R_l(t|\,x)\,,
$$
with
\begin{equation}
\cE^\hbar_{k,l}(x\,):=\,\int_{\R}\cE^\hbar(x,y)\, \Phi_k(y)\Phi_l(y)\,\dD y\,,\qquad k,\,l\,\in\,\N\,,
\label{def:Ekl}
\end{equation}
where the function $\cE^\hbar(.,.)$ is defined in \eqref{def:Eh}.

Let us first review some properties of $\cE_{k,l}$. Obviously $\cE^\hbar_{k,l}(x)\in\R$ and
\begin{equation*}
  % \label{Ekl:sym}
\cE_{k,l}(x) \,=\, \cE_{l,k}(x)\,,\qquad x\in\R, \,\,k,l\geq 0\,. 
\end{equation*}
Furthermore, applying the change of  variable $y\mapsto -y$ in
\eqref{def:Ekl} and  using that  $H_k(y)\,=\,(-1)^k H_k(-y)$ we also get
\begin{equation}
 \label{toto:1}
\cE^\hbar_{k,l}(x) \,=\, -(-1)^{k+l}  \,\cE^\hbar_{k,l}(x)\,,\qquad
x\in\R, \,k,l\geq 0\,. 
\end{equation}
This latter property implies that $\cE^\hbar_{k,l}$ is identically $0$
whenever $k+l$ is even.

We thus end up with the formulation of the von Neumann equation
\eqref{vNvW}  written in the Hermite basis in $y$ for all $k\geq 0$, 
  \begin{equation*}
%  \label{Hermite:D}
  \left\{
    \begin{array}{l}
  \ds\partial_t  R_{k} \,+\, i\,\left(\sqrt{\frac k 2}\,
      \cD\,R_{k-1}\,+ \, \sqrt{\frac{k +1}{2}}\,
      \cD^\star R_{k+1}\right)   \,=\,  -\,i\,\sum_{l\in\N} \cE^\hbar_{k,l}\,R_l\,,   
        \\[1.2em]
         \ds R_{k}(t=0) =  R^{in}_{\hbar,k}\,,
\end{array}\right.
\end{equation*}
where $\cD$ and $\cD^{\star}$ are given by
\begin{equation*}
%\label{AA}
\left\{
\begin{array}{ll}
\ds  \cD \,u   &=\,\ds +\partial_{x}  u \,+\, V^\prime(x)\, u\,, \\[1.1em]
\ds\cD^\star \,u  &=\,\ds -\partial_{x}  u \,+\, V^\prime(x)\, u\,.
\end{array}\right.
\end{equation*}
Let us emphasize an important property satisfied by $\cD$, which
we will need to recover later on, in the discrete setting. The
operator $\cD^\star$ is its adjoint operator in $L^2(\R)$, so that for
all $u$, $v\in {\rm Dom}(\cD)$ it holds
 \be 
 \label{prop1:A}
\left\langle\cD^\star u ,\, v\right\rangle = \left\langle u,\,\cD v\right\rangle,
\ee
 where $\langle.,\,.\rangle$ denotes the $L^2(\R^d)$ inner product.
 Then we define an approximation $R^{N}=(R_k^{N})_{0\leq k\leq N}
\in \mP_{N}(\R)$ solution to the
following system obtained after neglecting
 Hermite modes of order larger than $N$, that is,  for $0\leq k \leq {N}$ 
 \begin{equation}
  \label{Hermite:DN}
  \left\{
    \begin{array}{l}
  \ds\partial_t  R_{k}^{N} \,+\, i\,\left(\sqrt{\frac k 2}\,
      \cD R_{k-1}^{N}\,+ \, \sqrt{\frac{k +1}{2}}\,
      \cD^\star R_{k+1}^{N}\right)   \,=\,  -\,i\,\sum_{l= 0}^{N}\cE^\hbar_{k,l}\,R_l^{N}\,,   
        \\[1.2em]
         \ds R_{k}^{N}(0) =  R^{in}_{\hbar,k}\,,
\end{array}\right.
\end{equation}
where $R_{-1}^{N} = R_{N+1}^{N}=0$.

%%%%
%
% Modif
%
%%%% 

{Before to provide some error estimates,  we first
remark that the property \eqref{SAvW} now  becomes  for all $k\in\{0,\ldots, N\}$ and $t\geq 0$, 
  \begin{equation}
    \label{p:1}
  R^{N}_k(t) \,=\, (-1)^k\,\overline{R^{N}_k}(t)\,,
  \end{equation}
leading together with  \eqref{toto:1} on $\left( \cE_{k,l}^\hbar\right)_{k,l}$,  that for all $t\geq 0$, 
  \begin{equation}
    \label{p:2}
\| R^{N}(t) \|_{L^2} \,=\,  \| \cP_{N}R^{in} \|^2_{L^2}\,.
   \end{equation}
Furthermore, using again the property \eqref{toto:1},  the  trace property defined in \eqref{TRvW} now becomes 
$$
\int_\R   R^{N}(t|\,x,0)\,\dD x \,=\, \sum_{k=0}^{N}   \frac{H_k(0)}{\sqrt{2^k\,k!}}
                                                    \,  \int_\R 
                                   R^{N}_k(t|\,x) \,\dD x\,,
                                                    $$
but this quantity is not necessarily preserved with time $t\geq 0$ for the
truncated approximation.}

   \subsection{Review of error estimates}
\label{sec:2.2}
Here we suppose that for a fixed $m\geq
0$ and for all integer  $a$, $b$,
$\alpha$ and $\beta$ such that  $a+b+\alpha+\beta\le m$, we have 
\begin{equation}
\label{hyp:R}
\| x^a\partial^b_xy^\alpha\partial^\beta_y R^{in}_\hbar \|_{L^2}  \,<
\, \infty.
\end{equation}

In \cite{FG2:24}, we  have  obtained error estimates  between the
Hermite-Galerkin approximation and the solution to the von Neumann equation
\eqref{vNvW}.  As we have seen before, the formulation \eqref{Hermite:DN} is
well adapted  to obtaining  the semi-classical limit as $\hbar\rightarrow 0 .$ Indeed, the Hermite  formulation of the
semi-classical limit equation \eqref{eq:lim} is obtained simply
by neglecting the right-hand side term in
\eqref{Hermite:DN}. Thus,  $\hat R^N=(\hat R_k^N)_{0\leq k\leq  N}$ satisfies, for all $k\in\{0,\ldots,N\}$,
\begin{equation}
  \label{Hermite:limN}
    \left\{
    \begin{array}{l}
  \ds\partial_t  \hat R_{k}^N \,+\, i\,\left(\sqrt{\frac k 2}\,
      \cD\hat R_{k-1}^N\,+\, \sqrt{\frac{k +1}{2}}\,
      \cD^\star \hat R_{k+1}^N\right)   \,=\,  0\,,   
        \\[1.2em]
      \ds \hat R_{k}^N(t=0) =  \hat R^{in}_{\hbar, k}\,.
      
\end{array}\right.
\end{equation}
Moreover, under regularity assumption on the initial data, we prove error estimates of the Hermite-Galerkin
method  for smooth
solutions of the semi-classical limit equation \eqref{eq:lim}.

For a given $m\geq 0$ and for all integer  $a$, $b$,
$\alpha$ and $\beta$ such that  $a+b+\alpha+\beta\le m$, we suppose
that the initial data $\hat R^{in}_\hbar$ satisfies
\begin{equation}
\label{hyp:Rhat}
\| x^a\partial^b_xy^\alpha\partial^\beta_y \hat R^{in}_\hbar\|_{L^2}  \,<
\, \infty.
\end{equation}
Thus, the error estimates, proven in \cite{FG2:24}, can be summarized
in the following theorem.
\begin{theorem}
  \label{th:error}
  Let $p\geq 1$ and assume that $V$ satisfies
  \eqref{hyp:V1}-\eqref{hyp:V2} with $m=2(p+1)$ while  for all integers  $a$, $b$,
$\alpha$ and $\beta$ such that  $a+b+\alpha+\beta\le 2(p+1)$. Then we
have the following result :
\begin{itemize}
\item under the assumption \eqref{hyp:R} on the initial data,  the solutions $R$
to \eqref{vNvW} and  $R^{N}$ to \eqref{Hermite:DN} satisfy
\begin{eqnarray*}
\| R^N(t) - R(t) \|_{L^2} & \leq &   \frac{\cC_t}{(2{N}+3)^{p-1/2}}
                                   \,\left( 1 +\hbar^2 \right)\,,
\end{eqnarray*}
\item under the assumption  \eqref{hyp:Rhat}  on the initial data,  the solutions $\hat R$
  to \eqref{eq:lim} and  $\hat R^{N}$ to \eqref{Hermite:limN} satisfy
  \begin{eqnarray*}
&& \|\hat R(t)\,-\,\hat R^{N}(t)\|_{L^2} \,\leq\, \frac{\cC_t}{(2{N}+3)^{p-1/2}}\,,
  \end{eqnarray*}
  \end{itemize}
  where $\cC_t$ is a constant depending only on $p$, $V$ and the
initial data.
\end{theorem}

The key point in the proof of Theorem \ref{th:error} is to establish some  regularity
estimates on the solution to the von Neumann equation  written in
Weyl's variables

\section{Finite volume in $x$ and time discretization}
\label{sec:3}
\setcounter{equation}{0}
\setcounter{figure}{0}
\setcounter{table}{0}

This section is devoted to the discretization of the space
variable $x$
and  the time discretization. The aim is to preserve, at the discrete
level, the properties of the operators $\cD$ and $\cD^\star$ in \eqref{Hermite:DN}
and \eqref{Hermite:limN},  and the conservation of the $L^2$ norm for
both equations. Furthermore, in order to reduce the computational
cost, we also propose an approximation of the
non-local term $\cE^\hbar$ by only keeping only the second order term with respect to
$\hbar$. This latter approximation is valid only when $\hbar\ll 1$.

%%%%%%%%%%%%%%%%%%%%%%%%%%%%%%%%%%%%%%%%%%%%%%%%%%%%%%%%%%%%%%%%%%%%%%%%%%%
\subsection{Finite volume  method}
\label{sec:3.1}
We first treat the discretization with respect to the space
variable $x$,  providing an approximation of the operators $\cD$ and
$\cD^\star$. In order to get $L^2$ stability on the numerical solution, we
aim at preserving the duality structure \eqref{prop1:A}  of the discrete operators.

We fix a number of Hermite modes $N\in\N^*$ and consider an interval $(a,b)$ of $
\mathbb{R}$.  For $N_{x}\in\N^\star$, we introduce the set
$\cJ=\{1,\ldots, N_x\}$, and a family of control volumes
$\left(K_{j}\right)_{j\in\cJ}$ such that
$K_{j}=(x_{j-1/2},x_{j+1/2})$, where $x_j$ is the midpoint of the
interval $K_j$, {\it i.e.} $x_j=(x_{j-1/2},x_{j+1/2})/2$,  and 
\begin{equation*}
a=x_{{1}/{2}}<x_{1}<x_{{3}/{2}}<...<x_{j-{1}/{2}}<x_{j}<x_{j+{1}/{2}}<...<x_{N_{x}}<x_{N_{x}+{1}/{2}}=b\,.
\end{equation*}
Let us introduce  the mesh size
$$
\Delta x \,=\, x_{j+{1}/{2}}-x_{j-{1}/{2}}\,,
$$
and $\delta=(N,\Delta x)$ the numerical parameter. Then we  denote by 
$R_k^\delta(t)=(\cR_{k,j}(t))_{j\in\cJ}$, for $k\in\{0,\ldots,N\}$, the approximation of $R_k^N(t)$, where the index $k$
represents the $k$-th mode of the Hermite decomposition, whereas
$\cR_{k,j}(t)$ is an approximation of the mean value of $R_k(t)$
over the cell $K_{j}$ at time $t$.

First of all, the initial condition is discretized on each cell $K_{j}$ by:
\begin{equation*}
\cR_{k,j}({0})\,=\;\frac{1}{\Delta x}\int_{K_{j}}R_{k}^N(t=0,x)\, \dD x, \quad j\in\cJ\,.
\end{equation*}
By integrating equation
\eqref{Hermite:DN} on $K_{j}$ for $j\in\cJ$, we obtain the following
numerical scheme : 
\begin{equation}
	\frac{\dD R_{k}^\delta}{\dD t} \,+\, i \,\left(
	\sqrt{\frac{k}{2}}\,\cD_{\Delta x}\,R_{k-1}^\delta
        \,+\,\sqrt{\frac{k+1}{2}}\,\cD_{\Delta x}^\star\,R_{k+1}^\delta\right)
      \,=\,  -i\, \sum_{l=0}^N \cE^\hbar_{k,l}\, R_l^\delta \,,
	\label{discrete0}
      \end{equation}
 where $\cD_{\Delta x}$ (resp. $\cD_{\Delta x}^\star$)  is an
approximation of the operator $\cD$ (resp. $\cD^\star$) given by
\begin{equation*}
\cD_{\Delta x} = (\cD_j)_{j \in\cJ} \quad{\rm and}\quad \cD_{\Delta x}^\star = (\cD_j^\star)_{j \in\cJ}\, 
\end{equation*}
and where for
$R^\delta=(\cR_{j})_{j\in\cJ}$ it holds
\begin{equation}
\label{discrete1}
\left\{
  \begin{array}{l}
\ds\cD_{j} R^\delta \,=\, + \frac{\cR_{j+1} - \cR_{j-1}}{2\Delta
    x}  \,-\, E_j\, \cR_{j}\,, \quad j\in\cJ\,,
    \\[1.1em]
 \ds\cD_{j}^\star R^\delta \,=\, - \frac{\cR_{j+1} - \cR_{j-1}}{2\Delta
    x}  \,-\, E_j\, \cR_{j}\,, \quad j \in\cJ\,.
  \end{array}\right.
\end{equation}
On the other hand, the discrete force field $E_j$ is computed in terms
of the applied
potential $V$ as 
\begin{equation}
\label{discrete2}
E_j \,=\, -\frac{V(x_{j+1})-V(x_{j-1})}{2\Delta x}\,.
\end{equation}
It is worth mentioning that the same kind of approximation has been
successfully applied to the Vlasov-Poisson and
Vlasov-Poisson-Fokker-Planck systems. It allows to preserve the
structure of the continuous equation \cite{bf1-2024,bf2-2024, b-2024},
which is particularly convenient to capture the long time dynamics. 
Finally we apply  a Crank-Nicolson
scheme for the time discretization with a time step $\Delta t$.

\subsection{Approximation of the non-local term when $\hbar\ll 1$}
\label{sec:3.2}
Let us emphasize that the evaluation of the right-hand side
$\left(\cE^\hbar_{k,l}\right)_{0\leq k,\,l \leq N}$  in
\eqref{discrete0} may be costly since the computational complexity of
each time step is of order  $N^2\,N_x$. However, using the smoothness
of the potential $V$ and the fact that $\hbar\ll 1$, the non local
term can be drastically simplified.  Indeed, let us suppose that the potential  $V$ is real and analytic on
 $\R$, so that 
 $$
 V(z\,)=\,\sum_{n\ge 0}\frac1{n!}\,V^{(n)}(x)\,(z-x)^n\,,
 $$
hence we have
$$
\frac{V(x+\frac12\hbar y)-V(x-\frac12\hbar y)}{\hbar}\,=\,\sum_{n\ge 0}\frac{\hbar^{2n}}{2^{2n}\,(2n+1)!\,}V^{(2n+1)}(x)\,y^{2n+1}\,,
$$
and
\begin{equation*}
 % \label{DevTheta}
\cE ^\hbar_{k,l}(x)=\sum_{n\ge
  1}\frac{\hbar^{2n}}{2^{2n}(2n+1)!}V^{(2n+1)}(x)\, \int_\R
y^{2n+1}\,\Phi_k(y)\,\Phi_l(y)\,\dD y\,.
\end{equation*}
Applying twice the recursion relation \eqref{yPhi} shows that
$$
y^2\Phi_k(y)\,=\,\frac{\sqrt{k(k-1)}}2\,\Phi_{k-2}(y)\,+\,\frac{2k+1}2\,\Phi_k(y)\,+\,\frac{\sqrt{(k+1)(k+2)}}2\,\Phi_{k+2}(y)\,,
$$
Using again \eqref{yPhi} and the $L^2$-orthogonality of the Hermite
functions, we find that 
\begin{eqnarray*}
\int_\R
y^{3}\,\Phi_k(y)\,\Phi_l(y)\,\dD y&=&
                                      \sqrt{\frac{k(k-1)^2}8}\delta_{k,l+1}\,+\,\frac{2k+1}2\sqrt{\frac{k+1}2}\delta_{k,l-1}
  \\
                                  &+&\sqrt{\frac{(k+1)(k+2)(k+3)}8}\delta_{k,l-3}\,+\,\sqrt{\frac{k(k-1)(k-2)}8}\delta_{k,l+3}
  \\
  &+&\frac{2k+1}2\sqrt{\frac{k}2}\delta_{k,l+1}\,+\,\sqrt{\frac{(k+1)(k+2)^2}8}\delta_{k,l-1}\,.
\end{eqnarray*}
Combining  with the former Hermite-Galerkin method, we get
for $0\le k\le N$,
\begin{equation}
  \label{partie}
  \left\{
\begin{aligned}
& \ds\frac{\dD R_k^\delta}{\dD t} \,+\,
i\,\left( \sqrt{\frac{k}2}\cD_{\Delta x} R_{k-1}^\delta \,+\,
  \sqrt{\frac{k+1}2}\cD^\star_{\Delta x} R_{k+1}^\delta\right)
\,=\, -\,i\,\hbar^2 \,\cS_k[R^\delta]\,,
\\[1.1em]
&
\ds\cS_k[R^\delta] \,=\, \frac{V^{(3)} }{24}\,\left(  \frac32\,k\,\sqrt{\frac{k}2}\,R_{k-1}^\delta
\,+\, \frac32\,(k+1)\,\sqrt{\frac{k+1}2}\,R_{k+1}^\delta\right.\\
\\
& \qquad \qquad \left. +\, \sqrt{\frac{k(k-1)(k-2)}8}\,R_{k-3}^\delta\,+\,
  \sqrt{\frac{(k+1)(k+2)(k+3)}8}\,R_{k+3}^\delta \right),
\end{aligned}\right.
\end{equation}
where $ V^{(3)}(x)$ is an approximation of the third order derivative of
the potential $V$.

This approximation of the source term  introduces an additional error
of order $\hbar^4$, which means that this latter approximation is valid only when $\hbar$ is small
compared to $1/N$, where $N$ is the number of Hermite modes,  $\Delta x$ the mesh size and $\Delta t$ the time  time.

\section{Numerical simulations}
\label{sec:4}
\setcounter{equation}{0}
\setcounter{figure}{0}
\setcounter{table}{0}

In the following we perform
several numerical simulations to illustrate the efficiency of  the
proposed  Finite Volume/ Hermite Spectral method for the approximation
of the solution to the von Neumann equation \eqref{vNvW} for various
potentials $V$. We will for instance consider smooth potential with strong confining properties. 

\subsection{Harmonic potential}
\label{sec:4.1}

We first consider the  isotropic harmonic potential
$$
V(x) = \frac{x^2}{2}
$$
in order to benchmark the convergence rate of the proposed
method. Indeed,   the von Neumann equation \eqref{vNvW}  with this potential
has explicit solutions, which can be used as references to validate
numerical results. In this situation,  the Wigner equation \eqref{eq:Wigner} on
$W$ reduces to the transport equation in phase-space
$$
\partial_t W \,+\,\xi\partial_x W \,-\, \partial_x V\, \partial_\xi W\,=\,0\,.
$$
The exact solution is $W (t,x,\xi) = W_0(x(t), \xi(t))$, where
$$
\left\{
\begin{array}{l}
\ds x(t) \,=\, \cos(t) \,x(0)\,-\,  \sin(t)\,\xi(0)\,,
\\[0.9em]
\ds \xi(t) \,=\, \sin(t)\,x(0) \,+\, \cos(t)\,\xi(0)\,.
\end{array}\right.
$$
Consider a Gaussian wavepacket
$$
W^{in}(x,\xi) \,=\, \frac{1}{2\pi}\, \exp\left(-\frac{(x-5)^2 + \xi^2}{2}\right)
$$
as  initial condition, so that the wavepacket returns to the
initial state at the final time $T = 2\pi$. It is worth mentioning
that  our Hermite spectral method is well suited to compute the Wigner
transform.  Indeed, $W$ corresponds to the Fourier transform of $R$ and
Hermite functions $\Phi_k$ are eigenfunctions of the Fourier
transform : 
$$
\cF(\Phi_k)(\xi) \,=\, (-i)^k \,\Phi_k(\xi),  
$$
where $\cF$ represents the Fourier transform.

Here, we choose the $x$ domain as $[-8,8]$ and  mainly focus on the convergence with respect to $\Delta x$, $1/N$, and the time step $\Delta t$.  Numerical errors  are
presented in Table \ref{tab:1}, where the numerical error is given as
$$
\cE(\Delta t,\Delta x) \,=\, \max_{n\geq 0}\left(\sum_{k=0}^{N}
\sum_{j=0}^{N_x} \Delta x \, | \cR_{k,j}^n - R_k(t^n,x_j) |^2\right)^{1/2}. 
$$

From the results, we can make the following observations. On the one
hand, we present two tables showing the numerical error for $N=20$ and $N=200$. For
this choice of the initial data, few Hermite modes may be used and the
numerical error due to the truncation $N$ is negligible.  On the
other hand, the convergence rate
is plotted in the third column of each table.  One can achieve second order convergence in
$\Delta x$ and $\Delta t$, according with the theoretical value of the
centred finite volume scheme in $x$ and Crank-Nicolson time discretization. 

\begin{table}
  \begin{minipage}{0.5\linewidth}
    \centering
   	\begin{tabular}{|l|l|l|}
          \hline
	$(\Delta t,\,\Delta x) $  &  $L^2$ numerical error  & Order \\ \hline
	0.640 &    0.21458 & X  \\ \hline
	0.320  & 0.10280 & 1.06 \\ \hline
	 0.160  & 0.02930 & 1.82 \\ \hline
          0.080  &0.00766 & 1.94 \\ \hline
          0.040  &0.00191 & 2.00 \\ \hline		
	\end{tabular}
        \end{minipage}\hfill
        \begin{minipage}{0.5\linewidth}
          \centering
          \begin{tabular}{|l|l|l|}
          \hline
	 $(\Delta t, \,\Delta x) $
          &  $L^2$ numerical error & Order  \\ \hline
	0.640 &   0.22687 & X \\ \hline
	0.320  & 0.10379 &  1.13 \\ \hline
	 0.160  & 0.02931 & 1.81  \\ \hline
            0.080  & 0.00736 & 1.99 \\ \hline
            0.040  & 0.00184 &  2.00 \\ \hline
            \end{tabular}
          \end{minipage}
          \\[1.1em]
          \begin{minipage}{0.5\linewidth}
            \centering       
      (a) $N=20$
       \end{minipage}\hfill
       \begin{minipage}{0.5\linewidth}
         \centering
         (b) $N=200$ 
   \end{minipage}         
\caption{{\bf Harmonic potential :}  $L^2$ error with respect to $(\Delta t, \Delta x)$ and the number of Hermite modes $N$.}
 \label{tab:1}
 \end{table}

 %Finally in Figure  \ref{fig:11}, we present several snapshots of the
 %Wigner function in phase space $x-\xi$ with $N=20$ and $\Delta x=0.04$. 

 %\begin{figure}[!htbp]
%\begin{center}
 % \begin{tabular}{ccc}
 %  \includegraphics[width=5.2cm]{t1-t0.png} &
 %   \includegraphics[width=5.2cm]{t1-t3.png}&
 %   \includegraphics[width=5.2cm]{t1-t6.png} 
 %    \\
 %    $t=00$ &  $t =04.18$  & $t= 08.38$
 %   \\
 %   \includegraphics[width=5.2cm]{t1-t9.png} &
 %   \includegraphics[width=5.2cm]{t1-t12.png}&
 %   \includegraphics[width=5.2cm]{t1-t15.png}
 %   \\
  %   $t=12.54$ &  $t =16.76$  & $t= 20.88$ 
 % \end{tabular}
 % \end{center}
 % \caption{{\bf Harmonic potential} : Time evolution of the Wigner  function computed from the approximation of $R$. The lines represent the level set of the Hamiltonian $H(x,\xi) = (x^2 + \xi^2)/2$.}
%\label{fig:11}
%\end{figure}

%%%%%%%%%%%%%%%%%%%%%%%%%%%%%%%%%%%%%%%%%%%%%%%%%%%%%%%%%%%%%%%%%%%%%%%%%%%

\subsection{Quantum revivals}
\label{sec:4.2}

%%%%%%%%%%%%%%%%%%%%%%%%%%%%%%%%%%%%%%%%%%%%%%%%%%%%%%%%%%%%%%%%%%%%%%%%%%%

The quantum revival is a fascinating phenomenon in quantum mechanics
that involves periodic recurrences of the quantum wave function to its
original form during time evolution. These recurrences may occur
multiple times in space, leading to the formation of fractional
revivals, or they may occur almost exactly to the wave function's
original state at the beginning, also known as full revivals
\cite{manfredi1996theory, suh1991}.

To observe this phenomenon, we now use  the following quartic potential:
$$
V(x) = \frac{x^2}{2}  \,+\, \frac{\beta\,  x^4}{4},
$$
with $\beta = 0.5$. By injecting this potential in the  von Neumann equation \eqref{vNvW}, one
can easily verify that the von Neumann equation takes the simple form:
$$
\partial_t R \,=\,i \left(\partial_{xy} R \,-\,  V^\prime(x) \,y R\,-\, \frac{\beta \,\hbar^2}{4} \,x\,y^3 R\right)\,
$$
since derivatives of order larger than  five of the potential  vanish
identically.  We then choose the initial datum
\cite{manfredi1996theory}, {which is a coherent state (in appropriate coordinates)}
$$
R^{in}(x,y) \,=\,\frac{1}{\sqrt{2\pi}\sigma_x} \exp\left( -
 \frac{1}{2} \left(\frac{x^2}{\sigma_x^2} + y^2 \right) \right)\,,
$$
with $\sigma_x=0.6$, which corresponds to the following Wigner distribution
$$
W^{in}(x,\xi) \,=\,\frac{1}{2\pi\sigma_x} \exp\left( -
 \frac{1}{2} \left(\frac{x^2}{\sigma_x^2} + \xi^2 \right) \right)\,.
$$

Here, we pick the $x$ domain to be $[-4,4]$ with $N_x=400$, 
$N=400$ and $\Delta t=0.01$ and use the scheme \eqref{partie}, since
the local approximation of $\cE^\hbar$ is exact for such a potential.

\begin{figure}[!htbp]
\begin{center}
  \begin{tabular}{cc}
    \includegraphics[width=7.8cm]{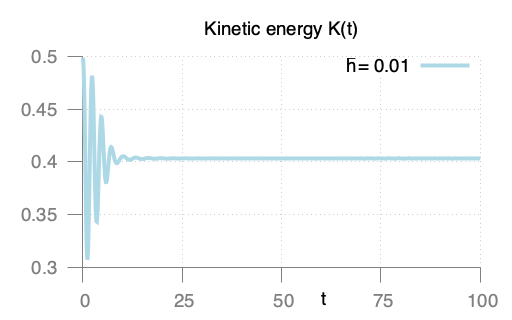} &
    \includegraphics[width=7.8cm]{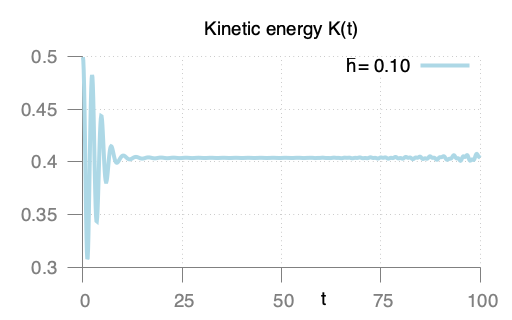}
\\
     \includegraphics[width=7.8cm]{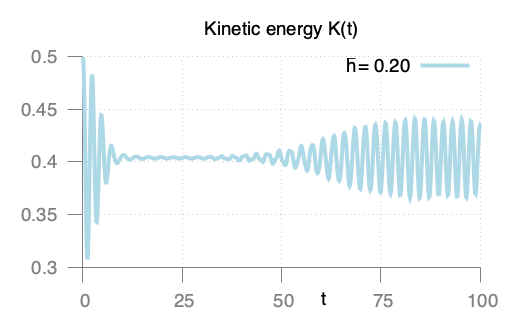} &
     \includegraphics[width=7.8cm]{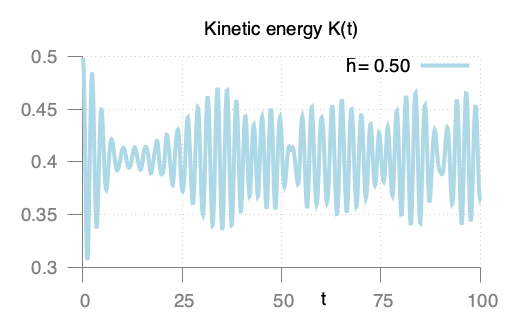}
   \end{tabular}
  \end{center}
  \caption{{\bf Quantum revivals} : Time evolution of the kinetic
    energy $\cK(t)$ for different values of $\hbar=0.01$, $\hbar=0.1$,
    $\hbar=0.2$ and $\hbar=0.5$.}
\label{fig:21}
\end{figure}

In \cite{manfredi1996theory}, quantum revivals, also named ``quantum
echoes'' have been studied in the context of an anharmonic
potential using the Wigner formalism, where quantum mechanics is
expressed in phase space. It was found  that mechanical effects prevent
complete phase mixing, leading to the emergence of a linear echo. Indeed, when $\hbar$ is small,  the kinetic energy relaxes to a
stationary value.  This is due to the phase mixing induced by the
quartic term in the potential $V$.  However, when $\hbar$ is large
($\hbar\geq 0.2$)  at a subsequent time, large oscillations (the
echo) appear. This phenomenon can be observed during the time
evolution of the kinetic energy,
$$
\cK(t) \,=\, -\frac{1}{2}\int_{\R} \partial_{y}^2 R(t,x,0)\,\dD x \,=\,
\int_{\R^2} W(t,x,\xi) \,  |\xi|^2\,\dD\xi\,\dD x\,, 
$$
which  is represented in Figure \ref{fig:21} for different values of
$\hbar=0.01$, $0.1$, $0.2$ and $0.5$. When $\hbar$
decreases, the echo amplitude goes to zero, while the time of its
appearance is rejected to infinity.  On the opposite, when
$\hbar=0.5$, the second order term  with respect to $\hbar$  in the
von Neumann equation becomes significant and prevents small structures
in phase space from emerging. As a result, the kinetic energy
maintains its oscillatory behavior  (see Figure \ref{fig:21}).

\begin{figure}[!htbp]
\begin{center}
  \begin{tabular}{ccc}
   \includegraphics[width=5.2cm]{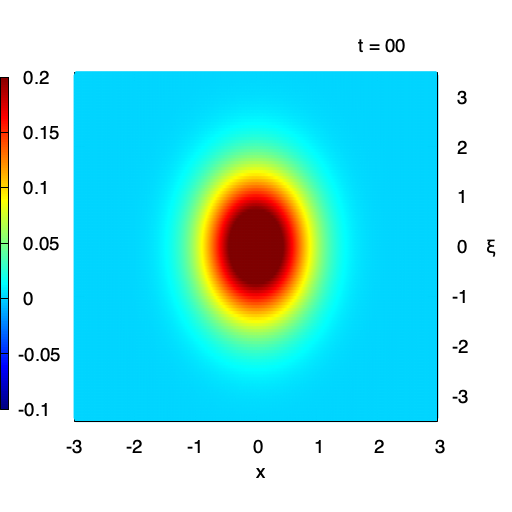} &
   \includegraphics[width=5.2cm]{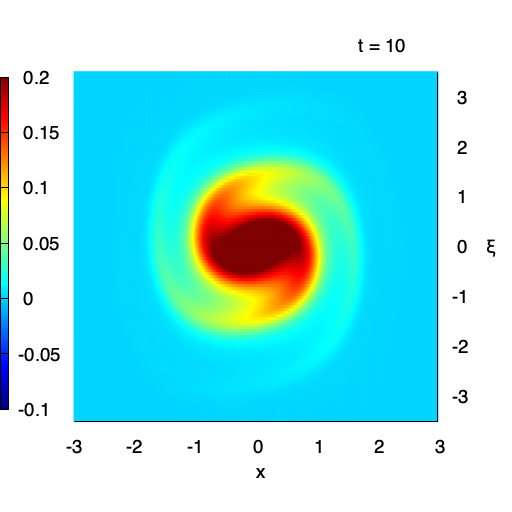} &
    \includegraphics[width=5.2cm]{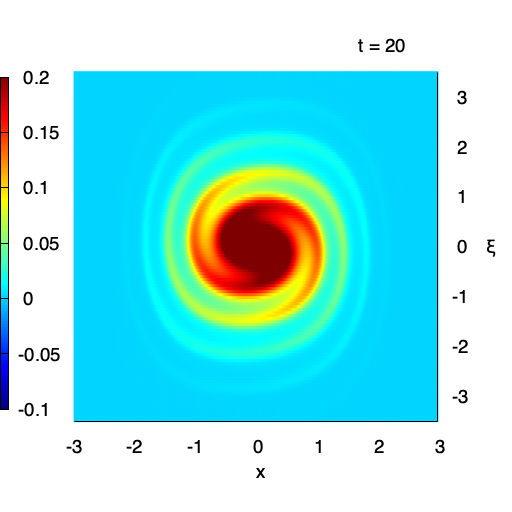}
    \\
   \includegraphics[width=5.2cm]{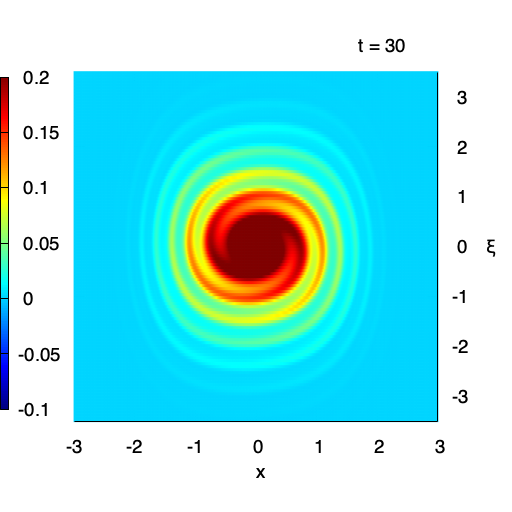} &
   \includegraphics[width=5.2cm]{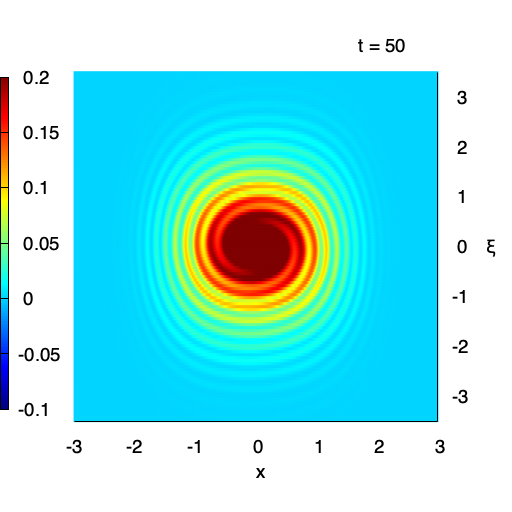} &
   \includegraphics[width=5.2cm]{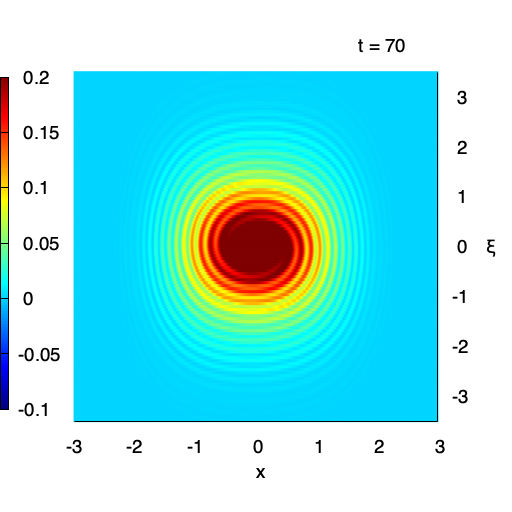}
   \end{tabular}
  \end{center}
  \caption{{\bf  Quantum revivals} : Snapshots of the Wigner
    distribution $W(t,x,\xi)$ at time  $t=0$, $10$, $20$, $30$, $50$
    and $70$ for $\hbar=0.01$.}
\label{fig:22}
\end{figure}

\begin{figure}[!htbp]
\begin{center}
  \begin{tabular}{ccc}
   \includegraphics[width=5.2cm]{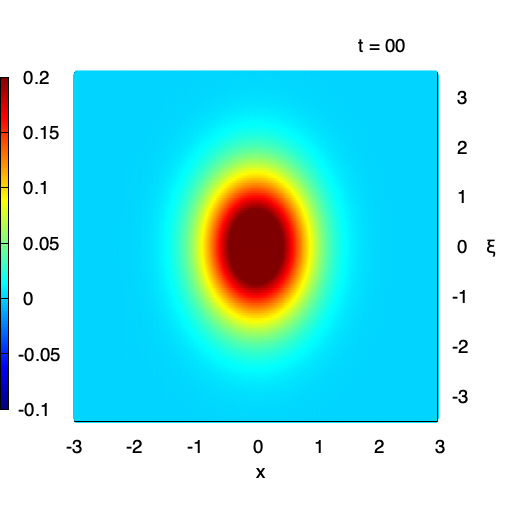} &
   \includegraphics[width=5.2cm]{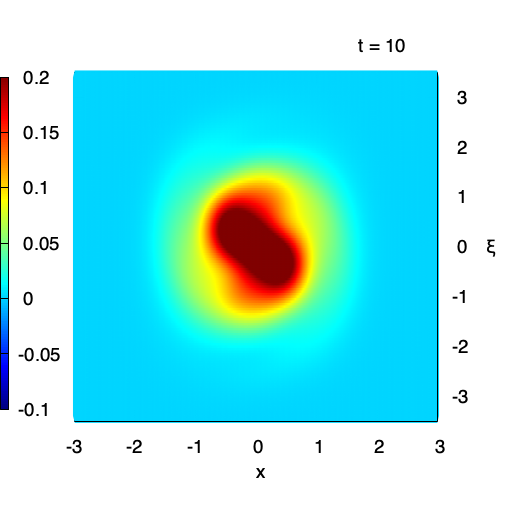} &
    \includegraphics[width=5.2cm]{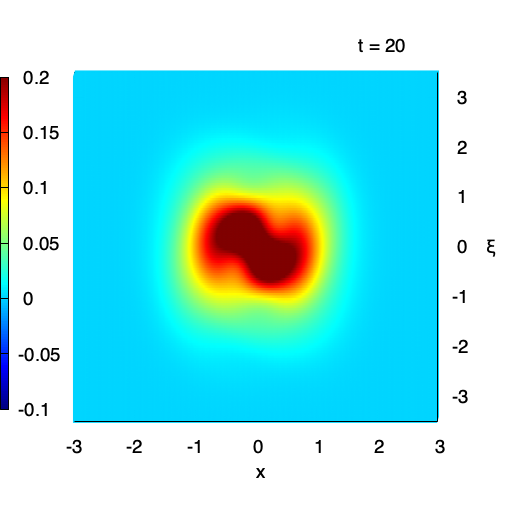}
    \\
   \includegraphics[width=5.2cm]{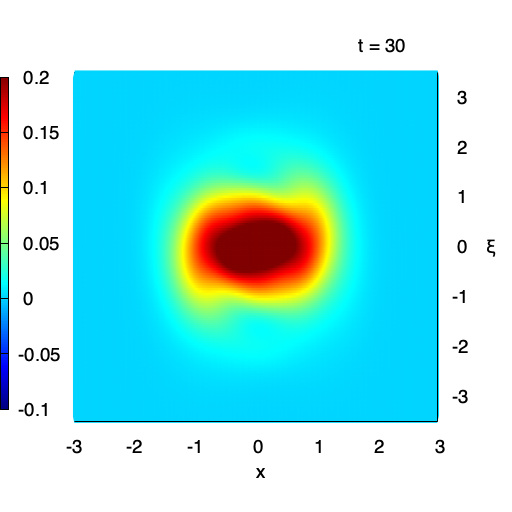} &
   \includegraphics[width=5.2cm]{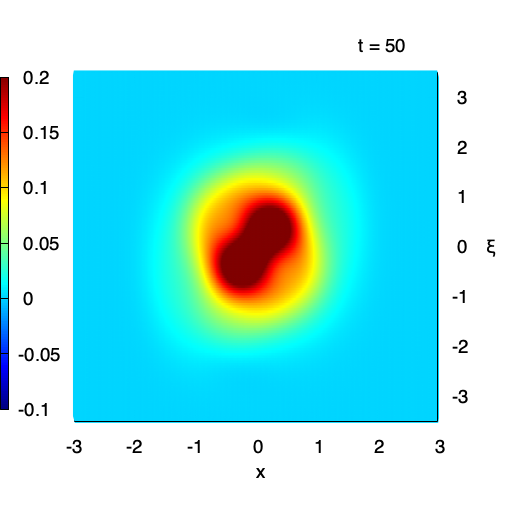} &
   \includegraphics[width=5.2cm]{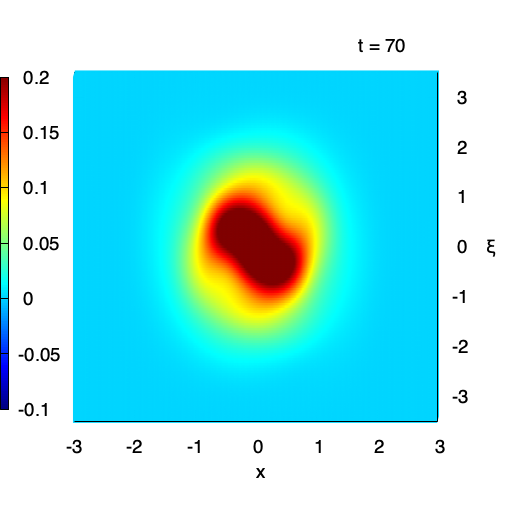}
   \end{tabular}
  \end{center}
  \caption{{\bf Quantum revivals} : Snapshots of the Wigner
    distribution $W(t,x,\xi)$ at time  $t=0$, $10$, $20$, $30$, $50$
    and $70$ for $\hbar=0.5$.}
\label{fig:23}
\end{figure}

Figures \ref{fig:22} and \ref{fig:23} show the phase portrait
for the evolution of the Wigner function starting from the  Gaussian initial state $W^{in}$. When $\hbar$
is small ($\hbar=0.01$), the distribution function $W$ soon develops a
spiral structure, displaying a very fine filamentation, which is the
ultimate cause of the kinetic energy relaxation. In the limit
$\hbar\rightarrow 0$ this phenomena still holds.  The same simulation is repeated in the quantum
case when $\hbar=0.5$ and is presented in  Figure \ref{fig:23}. We see
that complete phase mixing is now stopped by a nonzero
$\hbar$, via the second order term $\hbar^2 \,\beta\,x\,y^3 R/4$. The
correlations  among these structures are responsible for
the appearance of the echo. All these numerical results are consistent with
those presented in \cite{manfredi1996theory} on the numerical
simulation of the Wigner equation. It is worth mentioning
that  when $\hbar\ll 1$, the quartic potential predominates and  induces small
structures which are persistent in the limit $\hbar \rightarrow 0$,
hence  both  $\Delta x$ and $1/N$ have to be chosen  sufficiently small to follow the filamentation in
phase space (see Figure \ref{fig:22}) but this phenomena is not due to
the fact that $\hbar>0$!

\begin{figure}[!htbp]
\begin{center}
  \begin{tabular}{ccc}
   \includegraphics[width=5.2cm]{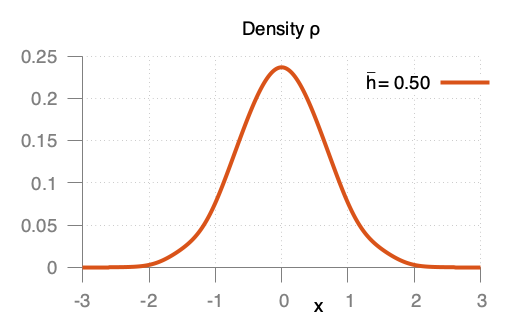} &
   \includegraphics[width=5.2cm]{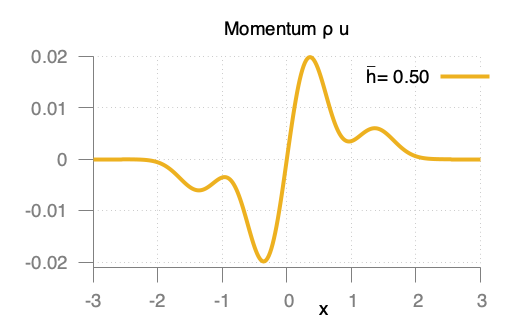} &
    \includegraphics[width=5.2cm]{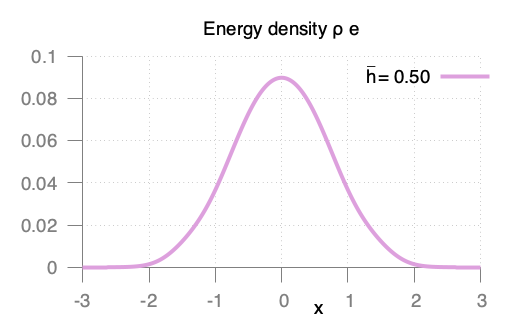}
    \\
    \includegraphics[width=5.2cm]{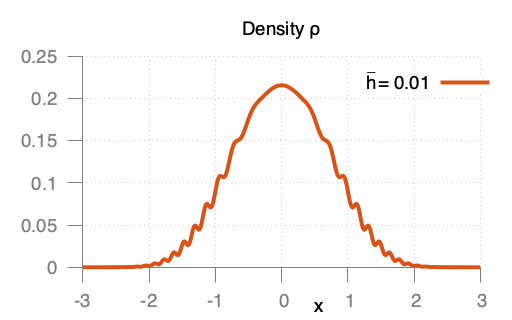} &
   \includegraphics[width=5.2cm]{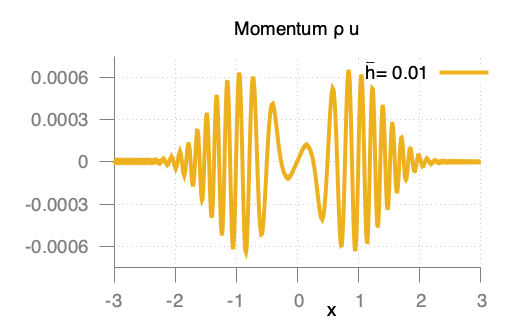} &
   \includegraphics[width=5.2cm]{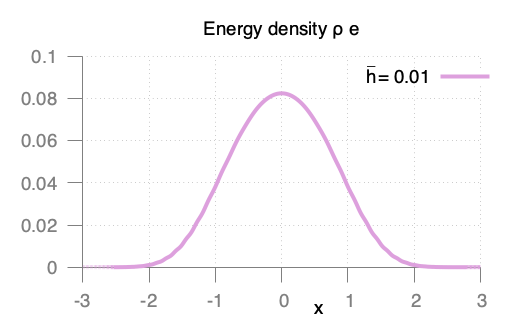}
       \\
    $\rho$ &  $\rho\, u$ & $\rho\, e$  
   \end{tabular}
  \end{center}
  \caption{{\bf Quantum revivals} : Snapshots of the matrix density
    $R(t,x,y)$ and the density $\rho$, momentum $\rho \;u$  and energy
    density $\rho\, e$ at time  $t=50$ for  $\hbar=0.5$ (top) and
    $\hbar=0.01$ (bottom).}
\label{fig:24}
\end{figure}

Finally in order to illustrate again the consequence of the amplitude
of $\hbar$ on the behavior of the solution, we present the macroscopic density $\rho$, momentum $\rho
u$ and $\rho e$, given by 
\begin{equation}
  \label{macro:q}
\left\{
  \begin{array}{l}
    \ds \rho(t,x) \,=\, R(t,x,0)\,,\\[0.9em]
    \ds \rho\,u(t,x) \,=\, i\,\partial_y R(t,x,0)\,,\\[0.9em]
    \ds \rho\,e(t,x) \,=\, -\frac{1}{2}\partial_y^2 R(t,x,0)\,,
    \end{array}\right.
\end{equation}
at time $t=50$ for $\hbar=0.5$ and $\hbar =0.01$ in Figure
\ref{fig:24}. When $0\leq \hbar \ll 0.01$, our numerical scheme is
stable and we do not observe the
influence of $\hbar$ on the solution since transport effects dominate.

While the energy density $\rho e$ remains consistent regardless of the
value of $\hbar$, the densities $\rho$ and $\rho u$ display distinct
behaviors. When $\hbar$ is small, the phase space filamentation
produces small oscillations in the density $\rho$ and the momentum
$\rho \,u$, while a significant oscillation is evident in the
momentum for $\hbar=0.5$.

\subsection{Quantum tunneling}
\label{sec:4.3}
Quantum tunneling is a phenomenon
in physics where particles like electrons or atoms can pass through
energy barriers that classical mechanics would deem impassable due to
insufficient energy. However, in quantum mechanics, there is a probability that a particle can ''tunnel'' through the barrier without having the required energy. Quantum tunneling is a crucial concept in various fields, including particle physics, solid-state physics, and nuclear fusion.

Here we choose the following potential $V$ :
$$
V(x) \,=\, e^{-x^2/2}\,,
$$
and  the initial datum in the von Neumann equation is the following {squeezed coherent state}
$$
R^{in}(x,y) \,=\,\frac{1}{\sqrt{2\pi}\sigma_x} \exp\left( -
 \frac{1}{2} \left(\frac{(x-x_0)^2}{\sigma_x^2} + y^2 \right)
\right)\, e^{4\,i\,y}\,,
$$
with $x_0=-5$ and $\sigma_x=0.6$.
Here, we pick the $x$ domain to be $[-12,16]$ with $N_x=1000$,
while  $N=400$ and $\Delta t=0.01$.

\begin{figure}[!htbp]
\begin{center}
  \begin{tabular}{ccc}
    \includegraphics[width=5.2cm]{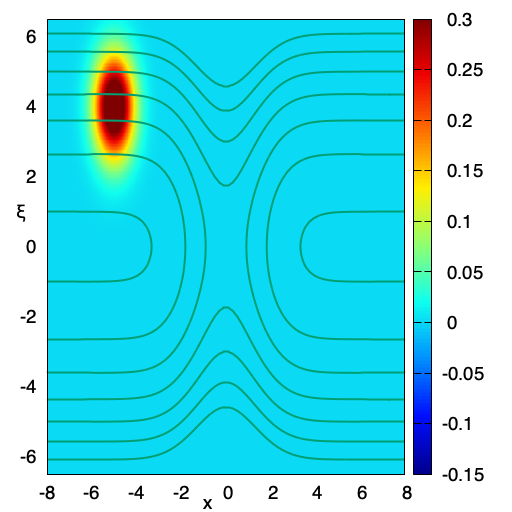} &
    \includegraphics[width=5.2cm]{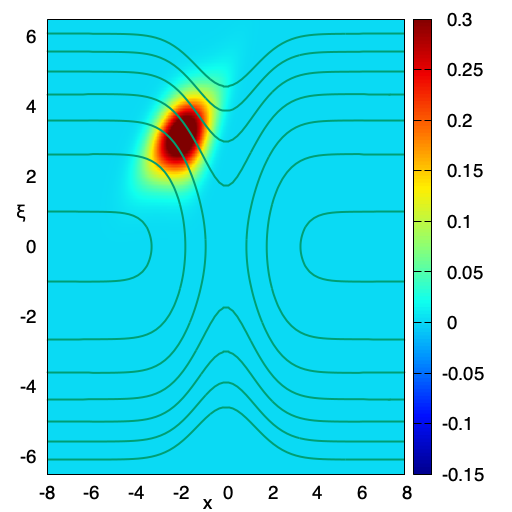}&
    \includegraphics[width=5.2cm]{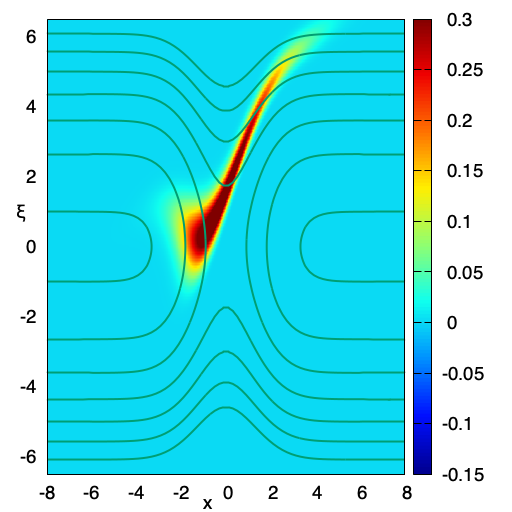} 
     \\
     $t = 0.0$ &  $t = 0.8$  & $t = 1.6$
    \\
    \includegraphics[width=5.2cm]{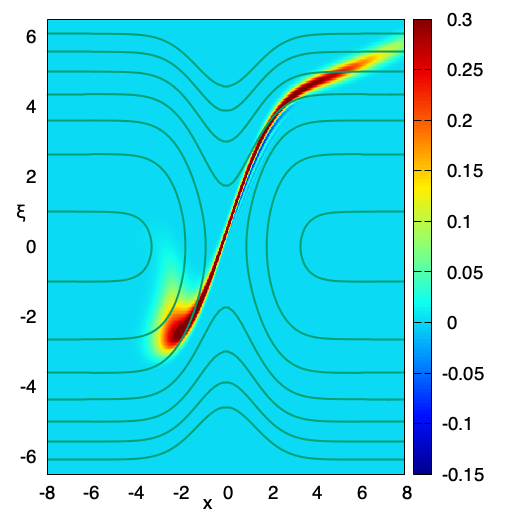} &
    \includegraphics[width=5.2cm]{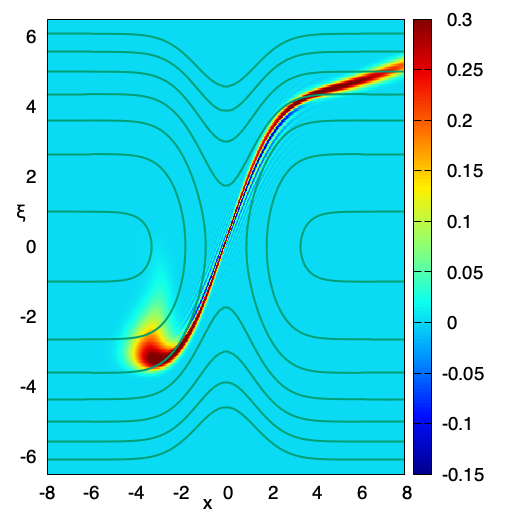}&
    \includegraphics[width=5.2cm]{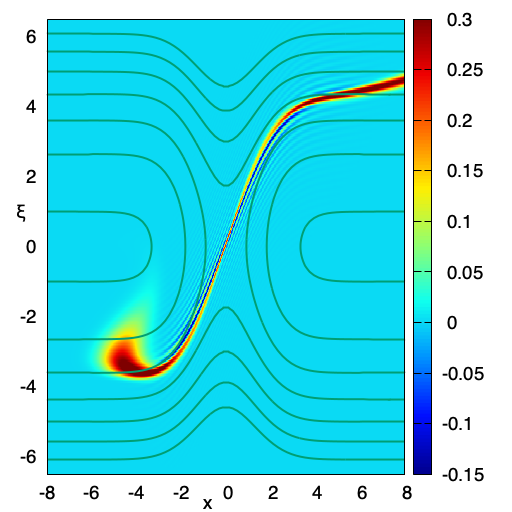}
    \\
     $t = 2.4$ &  $t =2.8$  & $t= 3.2$ 
  \end{tabular}
  \end{center}
  \caption{{\bf Quantum tunneling} : Time evolution of the Wigner
    function computed from the approximation of $R$ with $\hbar=0.1$. The lines
    represent the level set of the Hamiltonian $H(x,\xi) = V(x) + \xi^2/2$.}
\label{fig:31}
\end{figure}

\begin{figure}[!htbp]
\begin{center}
  \begin{tabular}{ccc}
   \includegraphics[width=5.cm]{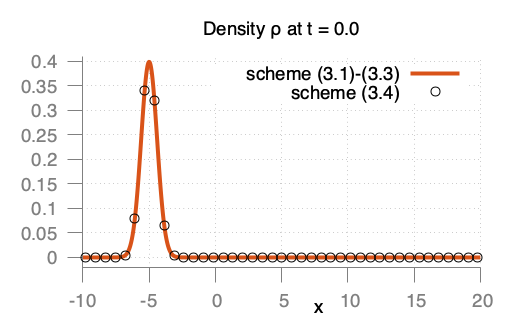} &
   \includegraphics[width=5.cm]{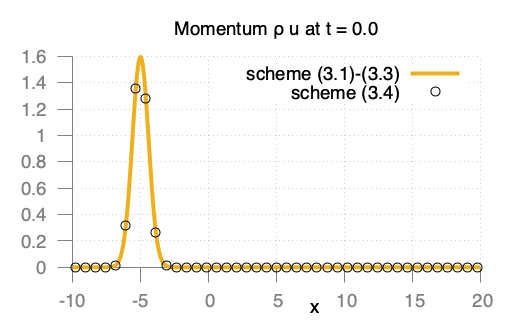} &
   \includegraphics[width=5.cm]{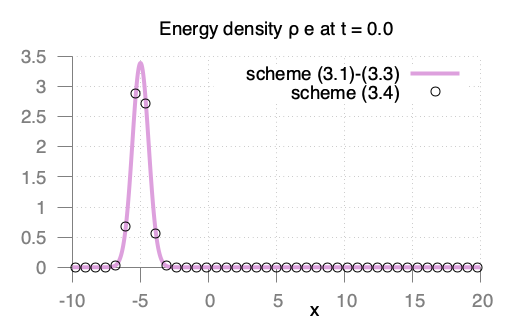}
    \\
    \includegraphics[width=5.cm]{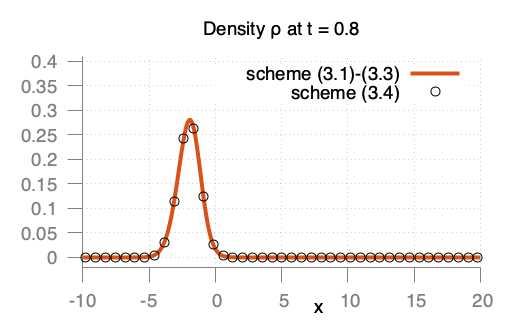} &
   \includegraphics[width=5.cm]{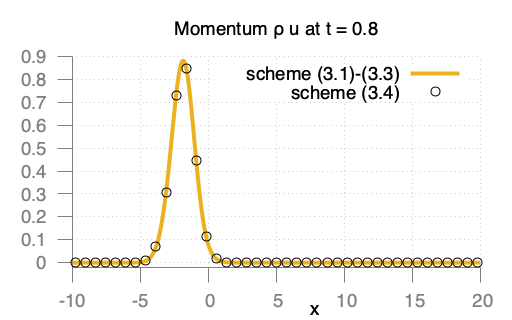} &
   \includegraphics[width=5.cm]{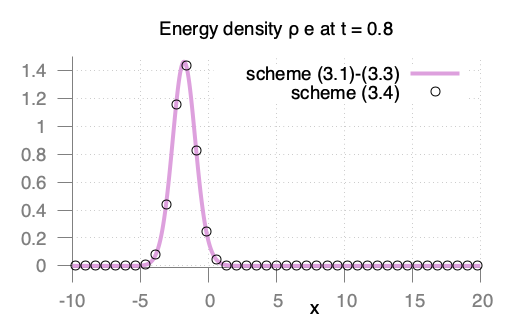}
    \\
   \includegraphics[width=5.cm]{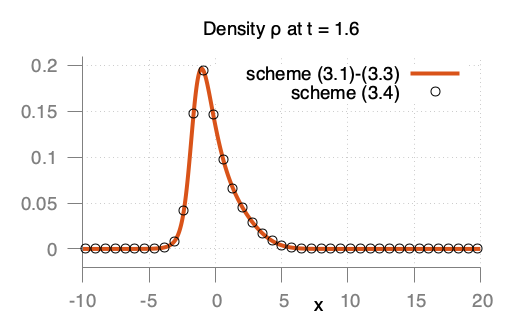} &
   \includegraphics[width=5.cm]{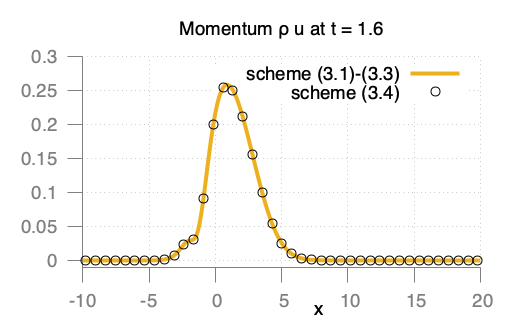} &
   \includegraphics[width=5.cm]{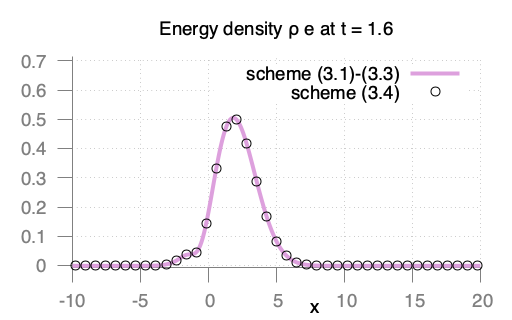}
    \\
   \includegraphics[width=5.cm]{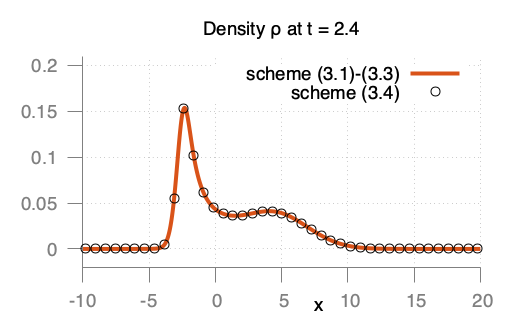} &
   \includegraphics[width=5.cm]{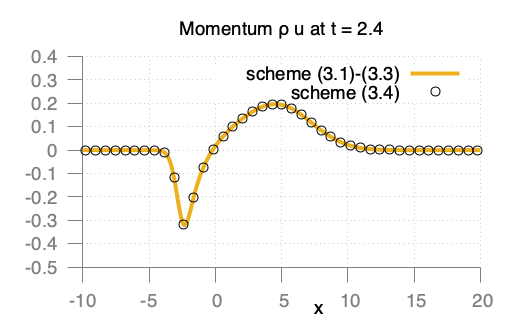} &
   \includegraphics[width=5.cm]{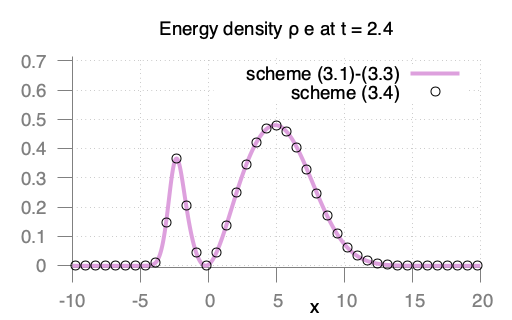}
    \\
   \includegraphics[width=5.cm]{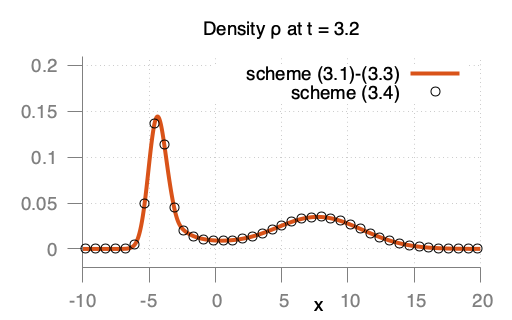} &
   \includegraphics[width=5.cm]{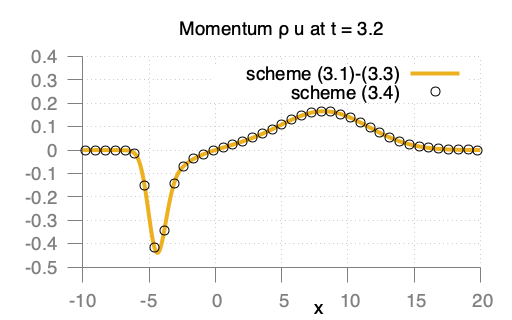} &
   \includegraphics[width=5.cm]{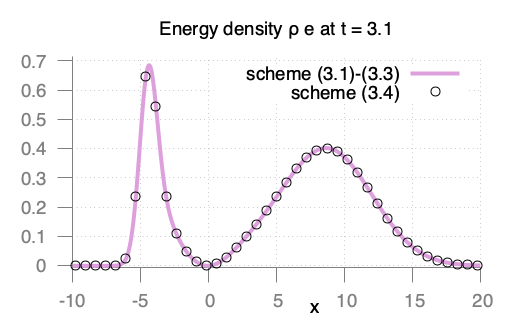}
   \\
    $\rho$ &  $\rho\, u$ & $\rho\, e$  
   \end{tabular}
  \end{center}
  \caption{{\bf Quantum tunneling} : Snapshots of the density $\rho$, momentum $\rho \;u$  and energy
    density $\rho\, e$ at time  different time. for $\hbar=0.1$}
\label{fig:32}
\end{figure}

On the one hand, we perform several numerical simulations with
$\hbar = 0.1$ and even smaller values ($\hbar=10^{-2}, \, 10^{-3}$) using
\eqref{discrete0}-\eqref{discrete2} where the nonlocal term
$\cE^\hbar$ is given by
\eqref{def:Ekl} and its truncated approximation \eqref{partie} of
order $\hbar^4$. For such a  configuration,  both approximations give similar
results and we only present those obtained with
\eqref{discrete0}-\eqref{discrete2}  and $\hbar=0.1$ in
Figures \ref{fig:31} and \ref{fig:32}.

The time evolution of the Wigner function in phase space, {\it i. e.}
in $(x,\xi)$-space, is proposed in Figure \ref{fig:31}. The initial Gaussian beam first moves forward, then it is deviated by the strong potential as it approaches $x=0$. Then at time
$t=1.6$, the beam goes backward but a small portion goes through the
barrier created by the potential and continues its way forward.
%%%%%%%%%%%%%%%%%%%%%%%%% 
%%%%%%%%%%%%%%%%%%
\begin{figure}[!htbp]
\begin{center}
  \begin{tabular}{ccc}
    \includegraphics[width=5.2cm]{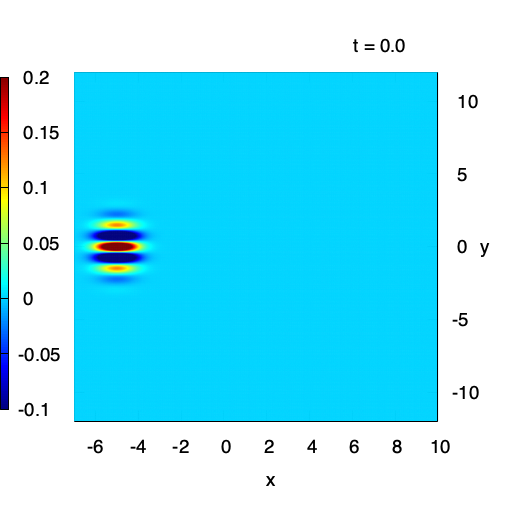} &
    \includegraphics[width=5.2cm]{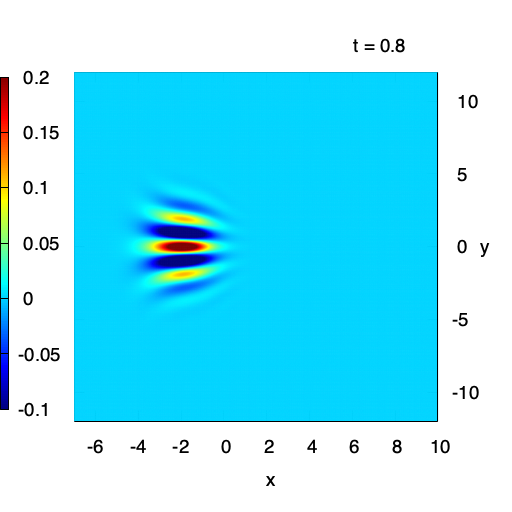}&
    \includegraphics[width=5.2cm]{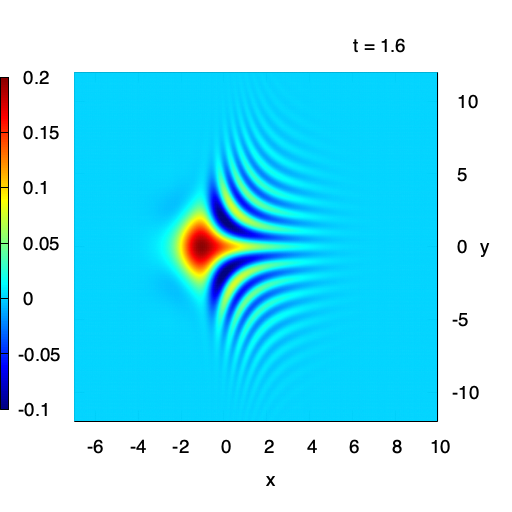} 
     \\
     $t = 0.0$ &  $t = 0.8$  & $t = 1.6$
    \\
    \includegraphics[width=5.2cm]{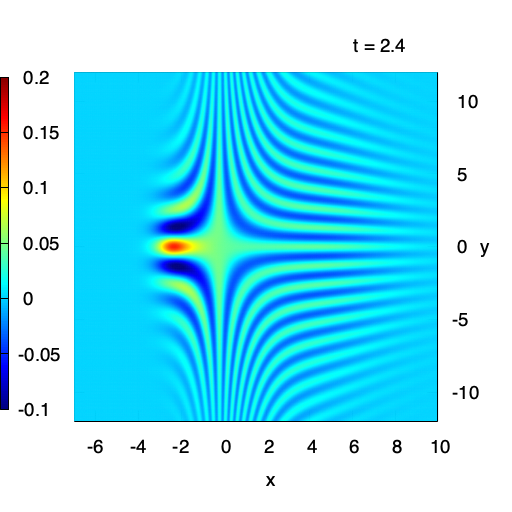} &
    \includegraphics[width=5.2cm]{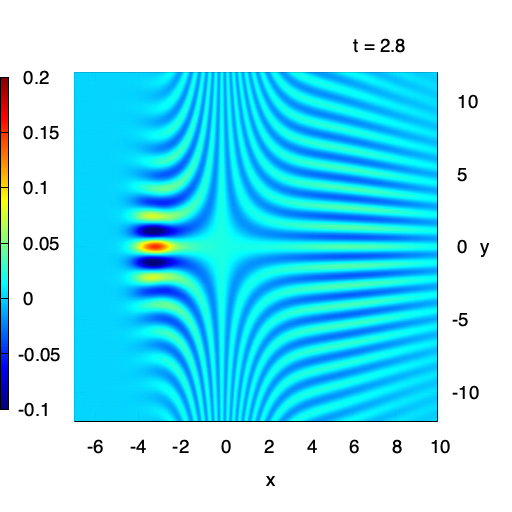}&
    \includegraphics[width=5.2cm]{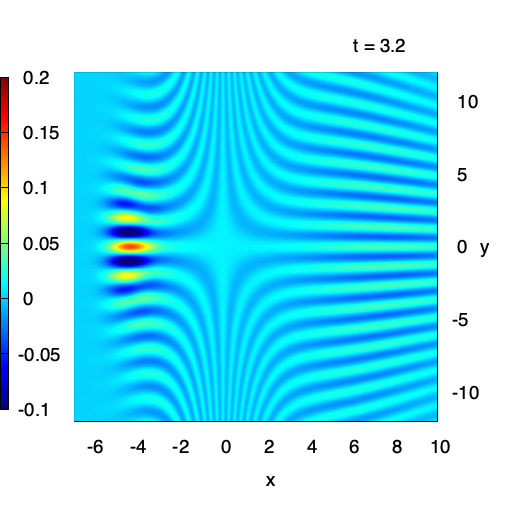}
    \\
     $t = 2.4$ &  $t =2.8$  & $t= 3.2$ 
  \end{tabular}
  \end{center}
  \caption{{\bf Quantum tunneling} : Time evolution of the approximation of Re$(R)$ with $\hbar=0.1$.}
\label{fig:31-bis}
\end{figure}
%%%%%%%%%%
\begin{figure}[!htbp]
\begin{center}
  \begin{tabular}{ccc}
    \includegraphics[width=5.2cm]{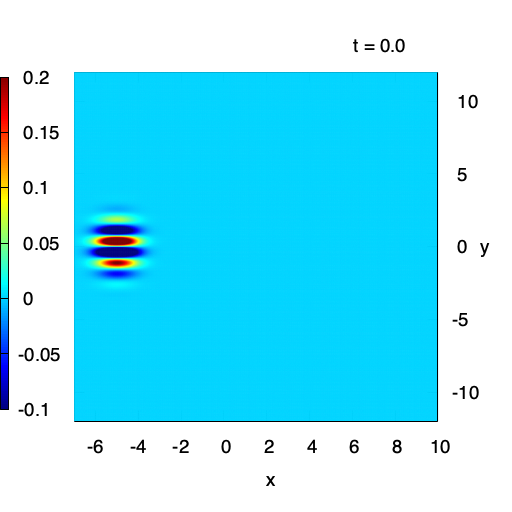} &
    \includegraphics[width=5.2cm]{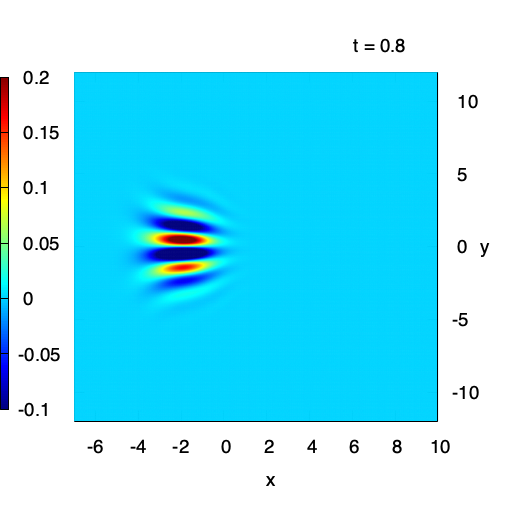}&
    \includegraphics[width=5.2cm]{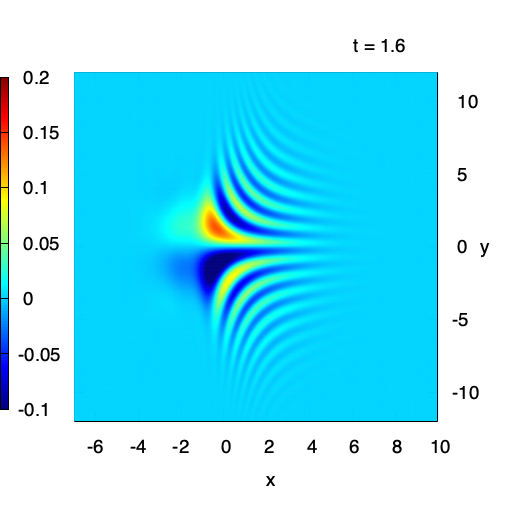} 
     \\
     $t = 0.0$ &  $t = 0.8$  & $t = 1.6$
    \\
    \includegraphics[width=5.2cm]{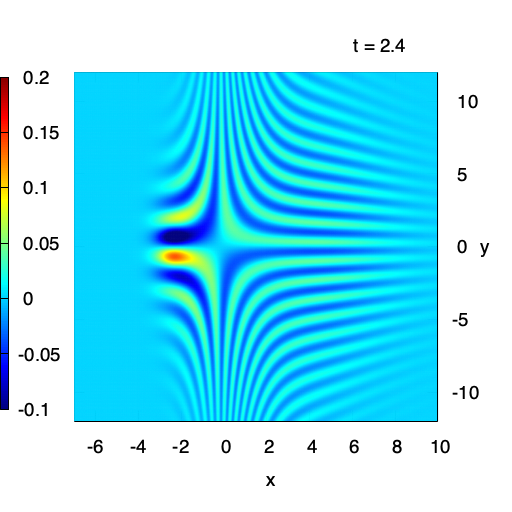} &
    \includegraphics[width=5.2cm]{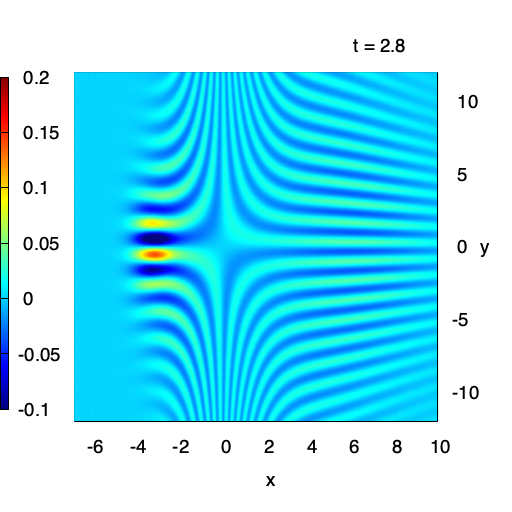}&
    \includegraphics[width=5.2cm]{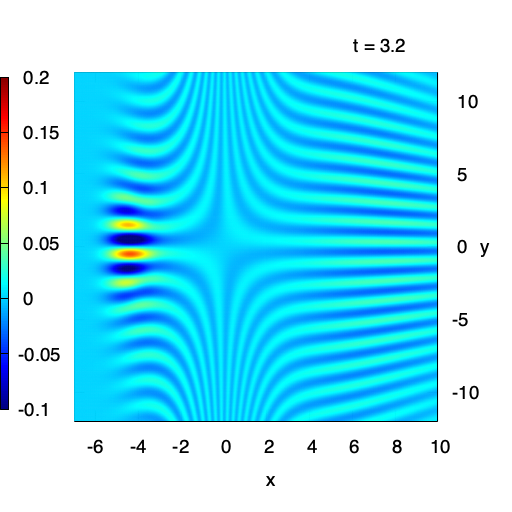}
    \\
     $t = 2.4$ &  $t =2.8$  & $t= 3.2$ 
  \end{tabular}
  \end{center}
  \caption{{\bf Quantum tunneling} : Time evolution of the approximation of Im$(R)$ with $\hbar=0.1$.}
\label{fig:31-ter}
\end{figure}
%%%%%%%%%%%%%%%%%%%%%%%%%%%%%%%%%%%%%%%%
This behavior can also be observed on the macroscopic quantities $\rho$, $\rho\, u$ and $\rho \,e$ defined in \eqref{macro:q} in Figure \ref{fig:32} when $\hbar =
0.1$. Indeed,  at time $t=1.6$ that a small amount of the density passes through
the quantum barrier and moves forward.
 When $\hbar\ll 0.1$, we do not observe any change on the numerical
approximation. Here again the use of Weyl's variables is a great advantage
in the vanishing $\hbar$  limit!

\begin{figure}[!htbp]
\begin{center}
  \begin{tabular}{cc}
    \includegraphics[width=7.8cm]{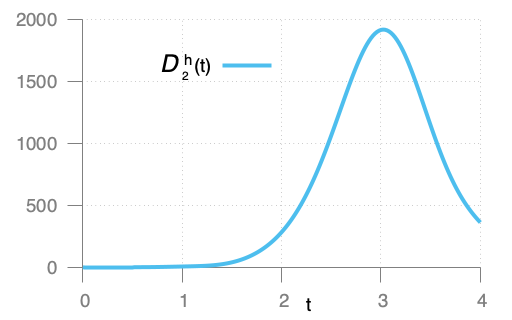} &
    \includegraphics[width=7.8cm]{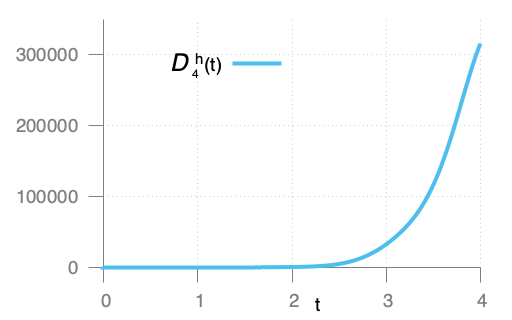}
    \\
    (a) $\cD_2^{\hbar}$ & (b) $\cD_4^{\hbar}$
   \end{tabular}
  \end{center}
  \caption{{\bf Quantum tunneling} : Time evolution of the quantities
    $\cD_2^{\hbar}$ and $\cD_4^{\hbar}$ for $\hbar=1$.}
\label{fig:33}
\end{figure}

On the other hand, when we increase the value of $\hbar \geq 0.5$, the local
approximation of  $\cE^\hbar$ given by
\eqref{def:Ekl} is not satisfying. To illustrate this fact we present in Figure \ref{fig:33},
the time evolution of the following quantities
{$$
\cD_{2}^\hbar(t) \, =\, \left\|\left(\frac{V\left(x+\frac\hbar 2\,
       y\right)\,-\,V\left(x-\frac\hbar 2\, y\right)}{\hbar}\,
        -\, V'(x) y\right)\,   R(t)\right\|_{L^2}
       $$}
and
       $$
\cD_{4}^\hbar(t) \, =\, \left\|\left(\frac{V\left(x+\frac\hbar 2\,
       y\right)\,-\,V\left(x-\frac\hbar 2\, y\right)}{\hbar}\,
       - \,\left(V'(x) y + \frac{\hbar^2}{24}\,V^{(3)}(x)\,
         y^3\right)\right)\, R(t)\right\|_{L^2} \,,   
$$
which are respectively  of order $\hbar^2$ and $\hbar^4$. Initially this
value is small but as the beam approaches the region where the
potential is strong, it increases rapidly and the local approximation
is not valid anymore.

\begin{figure}[!htbp]
\begin{center}
  \begin{tabular}{ccc}
    \includegraphics[width=5.2cm]{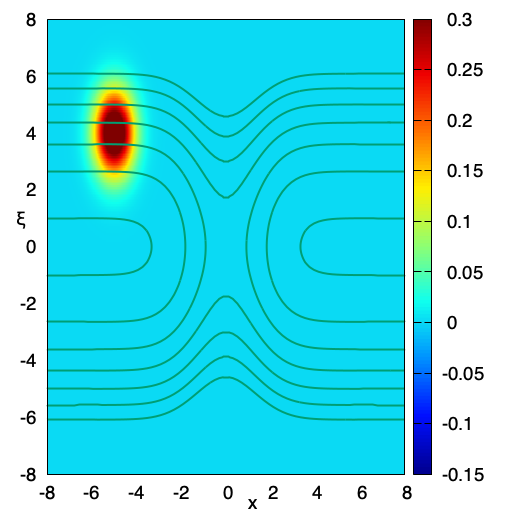} &
    \includegraphics[width=5.2cm]{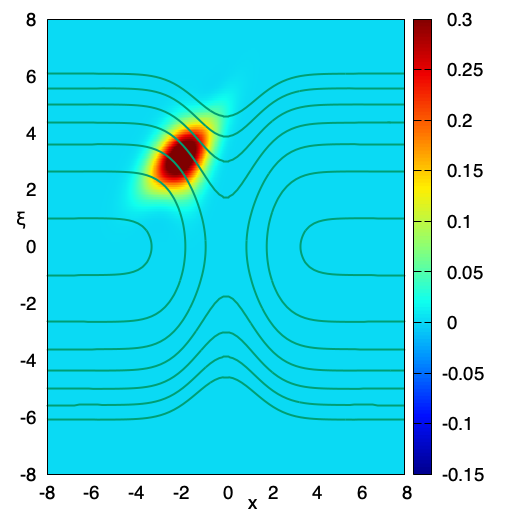}&
    \includegraphics[width=5.2cm]{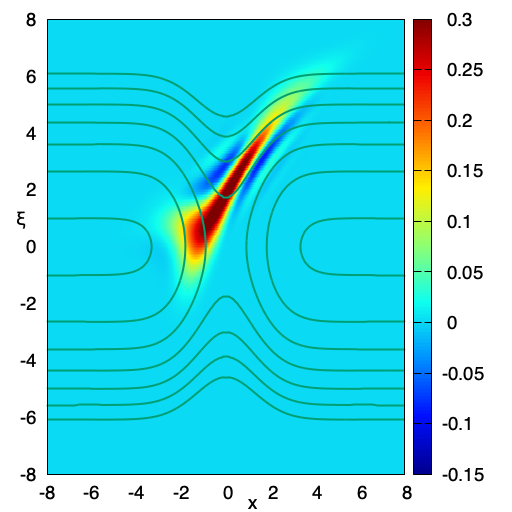} 
     \\
     $t = 0.0$ &  $t = 0.8$  & $t = 1.6$
    \\
    \includegraphics[width=5.2cm]{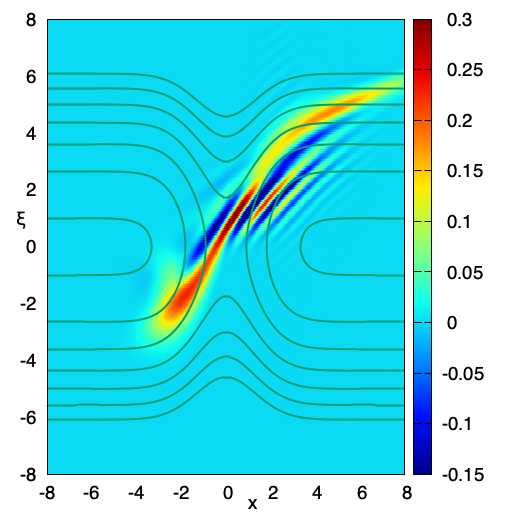} &
    \includegraphics[width=5.2cm]{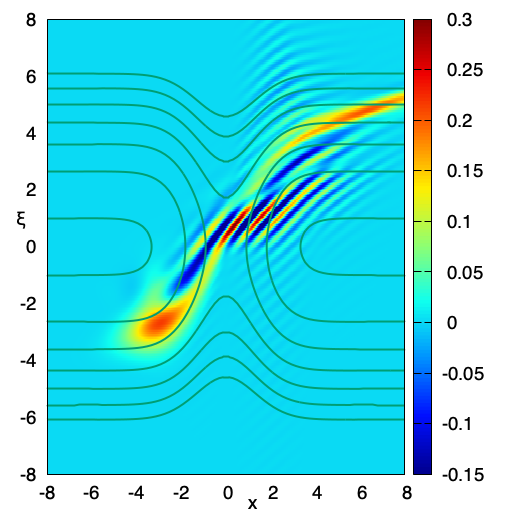}&
    \includegraphics[width=5.2cm]{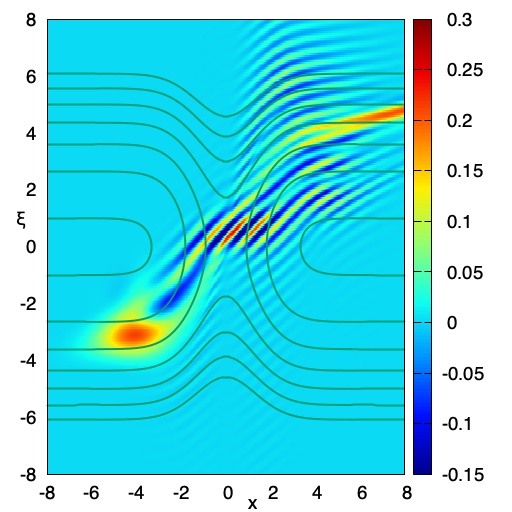}
    \\
     $t = 2.4$ &  $t =2.8$  & $t= 3.2$ 
  \end{tabular}
  \end{center}
  \caption{{\bf Quantum tunneling} : Time evolution of the Wigner
    function computed from the approximation of $R$ with $\hbar=1$. The lines
    represent the level set of the Hamiltonian $H(x,\xi) = V(x) + \xi^2/2$.}
\label{fig:34}
\end{figure}

Finally, we present the time evolution  of the Wigner function in
phase space, {\it i. e.}  $(x,\xi)$-space,  in Figure \ref{fig:34} for
$\hbar=1$. We see that it exhibits a  behavior slightly different from the one obtained with $\hbar=0.1$ since small waves appear
and propagate in the phase space domain.

\subsection{The Morse potential}
\label{sec:4.4}
\begin{figure}[!htbp]
\begin{center}
  \begin{tabular}{ccc}
    \includegraphics[width=5.2cm]{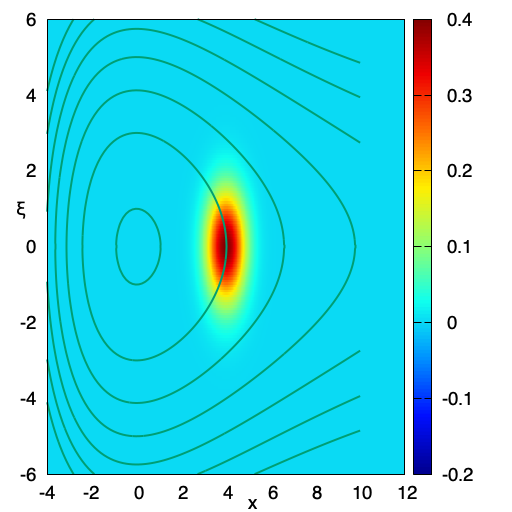} &
    \includegraphics[width=5.2cm]{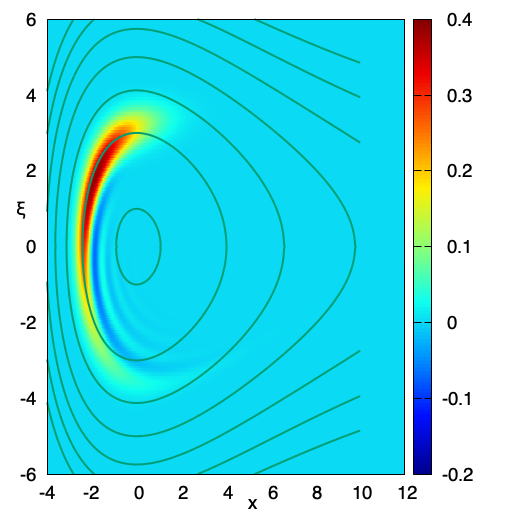}&
    \includegraphics[width=5.2cm]{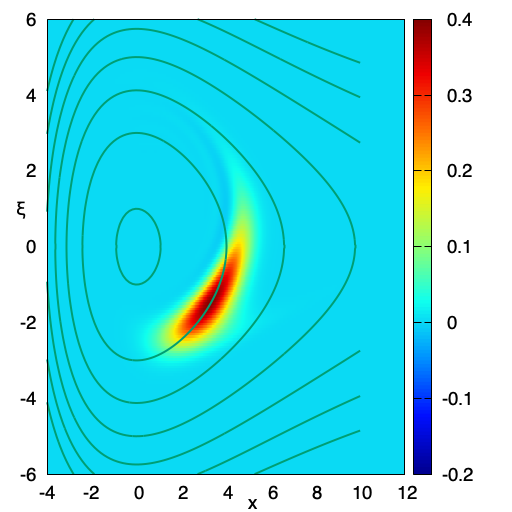} 
     \\
     $t = 00$ &  $t = 04$  & $t = 08$
    \\
    \includegraphics[width=5.2cm]{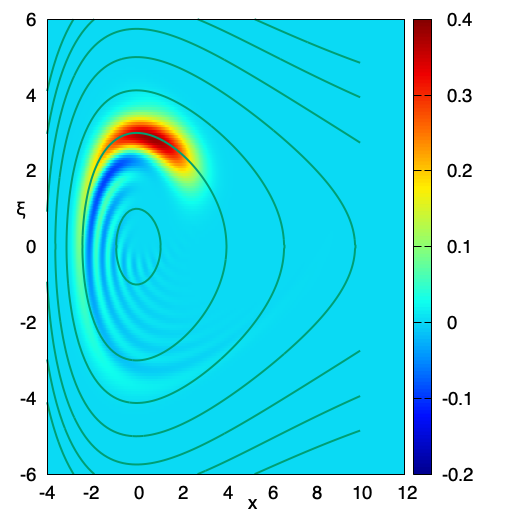} &
    \includegraphics[width=5.2cm]{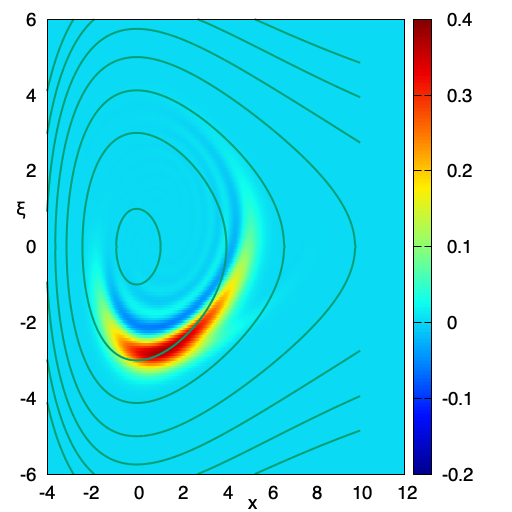}&
    \includegraphics[width=5.2cm]{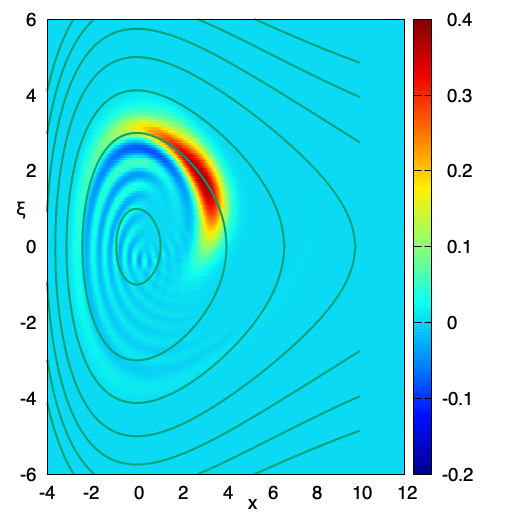}
    \\
     $t = 12$ &  $t =16$  & $t= 20$ 
  \end{tabular}
  \end{center}
  \caption{{\bf Morse potential} : Time evolution of the Wigner
    function computed from the approximation of $R$ with $\hbar=0.5$. The lines
    represent the level set of the Hamiltonian $H(x,\xi) = V(x) + \xi^2/2$.}
\label{fig:41}
\end{figure}

\begin{figure}[!htbp]
\begin{center}
  \begin{tabular}{ccc}
   \includegraphics[width=5.cm]{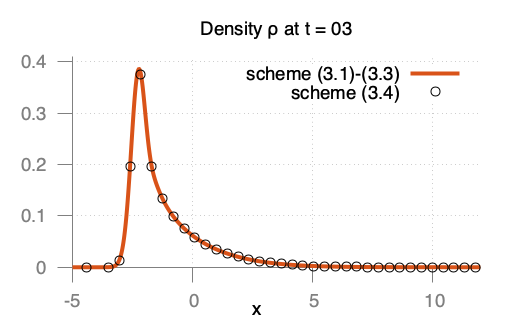} &
   \includegraphics[width=5.cm]{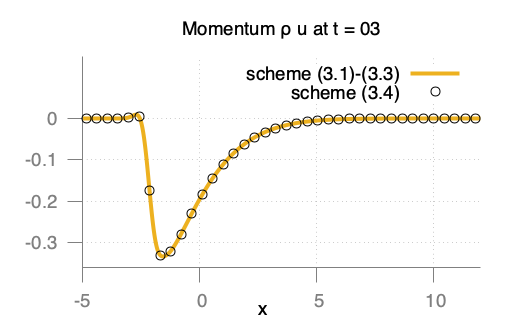} &
   \includegraphics[width=5.cm]{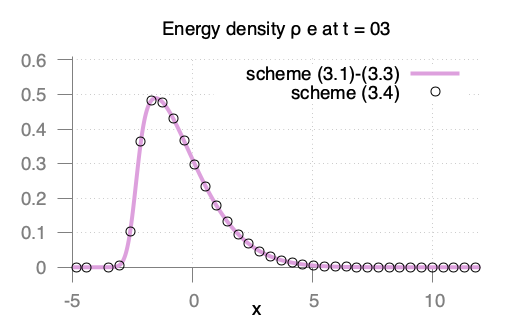}
    \\
    \includegraphics[width=5.cm]{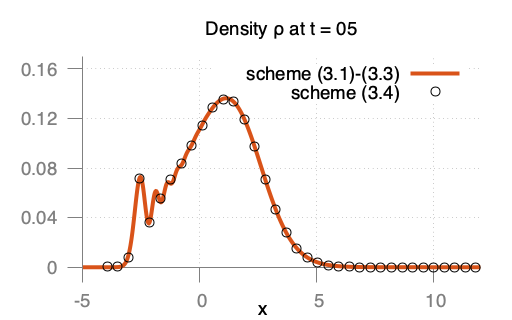} &
   \includegraphics[width=5.cm]{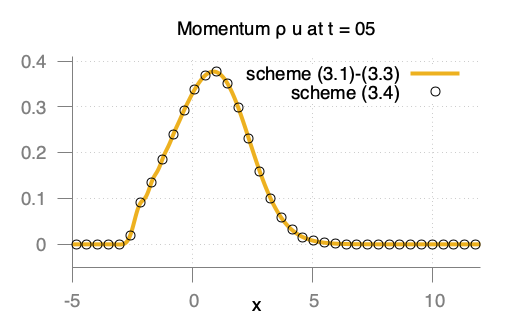} &
   \includegraphics[width=5.cm]{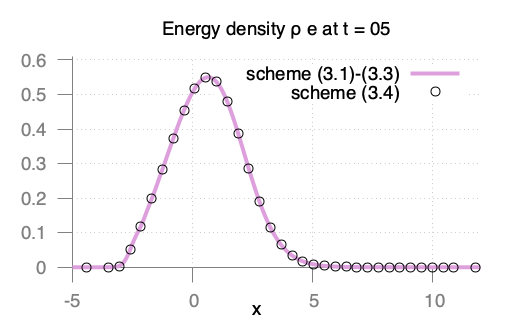}
    \\
   \includegraphics[width=5.cm]{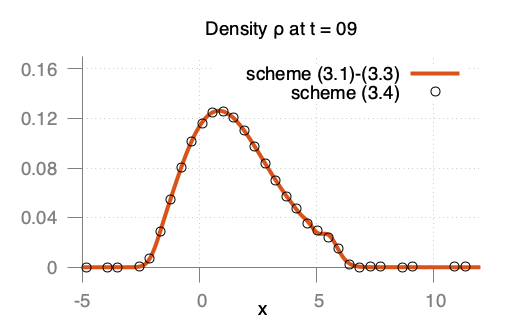} &
   \includegraphics[width=5.cm]{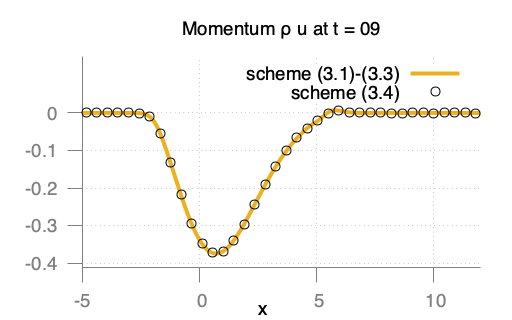} &
   \includegraphics[width=5.cm]{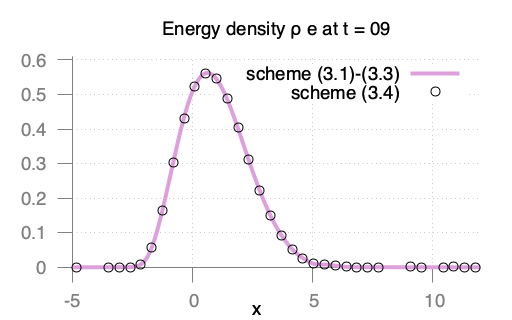}
    \\
   \includegraphics[width=5.cm]{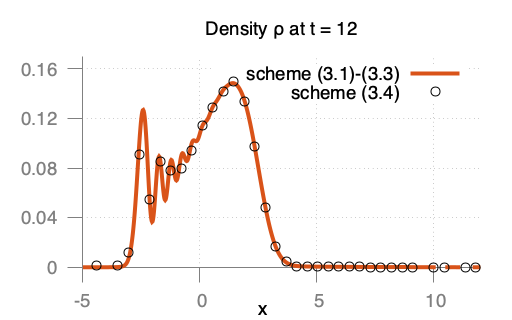} &
   \includegraphics[width=5.cm]{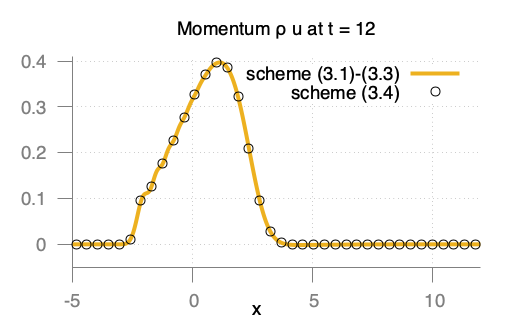} &
   \includegraphics[width=5.cm]{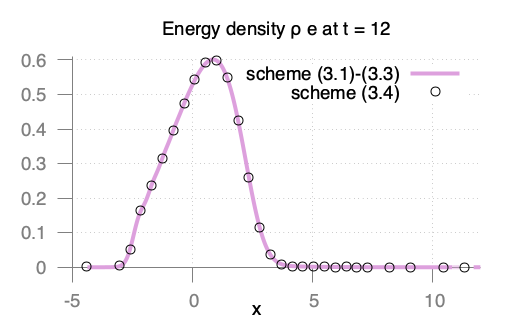}
    \\
   \includegraphics[width=5.cm]{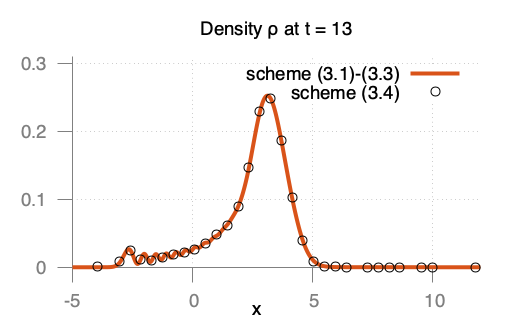} &
   \includegraphics[width=5.cm]{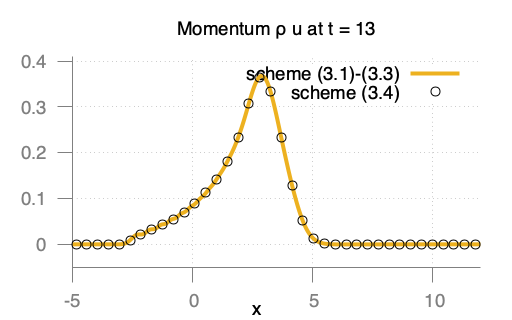} &
    \includegraphics[width=5.cm]{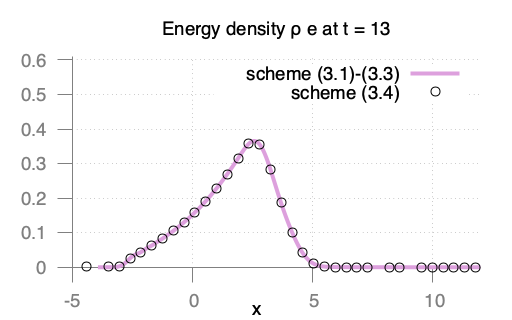}
    \\
   \includegraphics[width=5.cm]{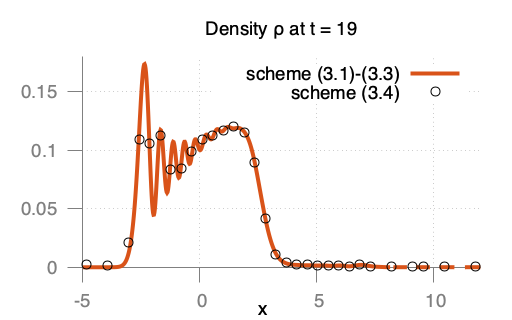} &
   \includegraphics[width=5.cm]{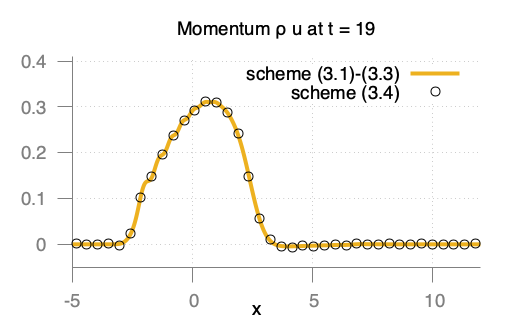} &
   \includegraphics[width=5.cm]{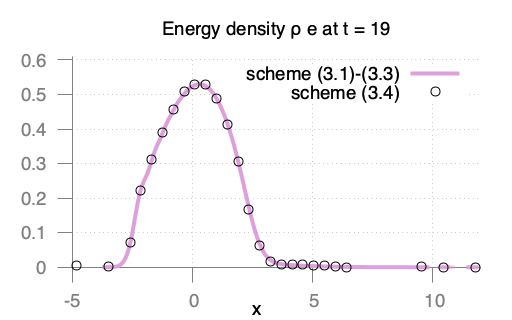}
   \\
    $\rho$ &  $\rho\, u$ & $\rho\, e$  
   \end{tabular}
  \end{center}
  \caption{{\bf Morse potential} : Snapshots of the density $\rho$, momentum $\rho \;u$  and energy
    density $\rho\, e$ at time  different time. for $\hbar=0.5$ at
    time $t=3$, $5$, $9$, $12$, $13$ and $19$.}
\label{fig:42}
\end{figure}

The Morse potential is used to approximate the vibrational structure
of a diatomic molecule, since it explicitly includes the effects of
bond breaking, such as the existence of unbound states. It also
accounts for the anharmonicity of real bonds and the non-zero
transition probability. The Morse potential energy function is of the
form
$$
V(x) = 20 \, \left( 1- \exp(-0.16 \,x)\right)^2\,,
$$
and  the initial datum in the von Neumann equation is {
  a coherent state (in appropriate coordinates)}
$$
R^{in}(x,y) \,=\,\frac{1}{\sqrt{2\pi}\sigma_x} \exp\left( -
 \frac{1}{2} \left(\frac{(x-x_0)^2}{\sigma_x^2} + y^2 \right)
\right)\,,
$$
with $x_0=4$ and $\sigma_x=0.6$.
Here,  the $x$ domain is $[-4,16]$ with $N_x=1000$, while
 $N=400$, $\Delta t=0.01$ and $\hbar=0.5$. We performed again  numerical
 computations using \eqref{discrete0}-\eqref{discrete2}, with nonlocal term $\cE^\hbar$ given by
\eqref{def:Ekl} and its truncated approximation \eqref{partie} of
order $\hbar^4$. Surprisingly both results are in good agreement for
the Wigner distribution and the macroscopic quantities (see Figures
\ref{fig:41} and \ref{fig:42}).

The time evolution of the Wigner function in phase space is proposed in Figure \ref{fig:41} for $\hbar=0.5$. The initial Gaussian beam
is rapidly deviated due
to the strong potential. We first see that the trajectory of the  beam
coincides with the isolines of the Hamiltonian
$$
H(x,\xi) \,=\, \frac{\xi^2}{2} \,+\, V(x)
$$
then the Gaussian beam is distorted and small waves appear in the
center and are
captured by the strong potential ($t=20$). These oscillatory waves can
be also observed on the density $\rho$ (see left column of Figure
\ref{fig:42}).

%%%%%%%%%%%%%%%%%%%%%%%%%%%%%%%%%%%%%%%%%%%%%%%%%%%%%%%%%%%%%%%%%%%%%%%%%%%

\section{Conclusion}
\label{sec:5}
In conclusion, we have introduced a new approach to discretize the
von Neumann equation in the semi-classical limit. By utilizing Weyl's
variables and employing a truncated Hermite expansion of the density
operator, we are able to address the stiffness associated
with the equation. Our method allows for error estimates by leveraging
the propagation of regularity on the exact solution. Therefore, our asymptotic
preserving numerical approximation, coupled with Hermite polynomials,
represents an interesting approach for accurately solving the von
Neumann equation in the semi-classical regime. Finally through the
development of a finite volume approximation and numerical
simulations showcasing phenomena in both quantum mechanics and its
semi-classical limit, we have shown the efficiency of the present
approach to describe mixed states with smooth Wigner
function. {The crucial point in our approach is the use
of Weyl's variable since it removes the stiffness due to the smallness
of $\hbar$, hence alternative methods based on stochastic algorithms \cite{XYZZ}
may be also applied in this context.}

\bibliographystyle{abbrv}
\bibliography{refer}

\end{document}